\documentclass[11pt]{article}

\newcounter{remarknum}
\makeatletter
\def\remark{\@ifnextchar[{\@with}{\@without}}
\def\@with[#1]{\refstepcounter{remarknum}\textit{Remark~\theremarknum.} \label{rem:#1}}
\def\@without{\refstepcounter{remarknum}\textit{Remark~\theremarknum.} }
\makeatother

\RequirePackage[a4paper, left=25mm, right=25mm, top=25mm, bottom=30mm]{geometry}
\RequirePackage[OT1]{fontenc}
\RequirePackage[british]{babel}
\RequirePackage{amssymb, amsmath, amsthm}
\allowdisplaybreaks[2]
\RequirePackage[sort&compress, numbers]{natbib} 

\DeclareRobustCommand{\NLprefix}[3]{#2}

\RequirePackage{hyperref}
\RequirePackage{enumerate}
\RequirePackage{color}
\RequirePackage[margin=1cm]{caption}
\RequirePackage{subcaption}
\RequirePackage{tikz}
\usepackage{tikz-qtree}
\usetikzlibrary{trees} 

\RequirePackage{calc}

\numberwithin{equation}{section}
\theoremstyle{plain}
\newtheorem{theorem}{Theorem}
\newtheorem{lemma}[theorem]{Lemma}
\newtheorem{corollary}[theorem]{Corollary}
\numberwithin{theorem}{section}

\DeclareSymbolFont{bbold}{U}{bbold}{m}{n}
\DeclareSymbolFontAlphabet{\mathbbold}{bbold}

\newcommand{\numberthis}{\stepcounter{equation}\tag{\theequation}}
\newcommand{\alert}[1]{\textcolor{red}{#1}} 

\renewcommand{\cite}{\alert{cite[?]}}
\renewcommand\citet\Citet
\renewcommand\citep\Citep
\renewcommand\citealt\Citealt
\renewcommand\citealp\Citealp
\renewcommand\citeauthor\Citeauthor

\newcommand{\levy}{L{\'e}vy }

\renewcommand{\a}{\alpha}

\renewcommand{\d}{\delta}
\newcommand{\dd}{\,\mathrm{d}}
\newcommand{\D}{\Delta}
\newcommand{\e}{\varepsilon}
\newcommand{\E}{\mathbb{E}}
\newcommand{\g}{\gamma}
\newcommand{\ind}{\mathbbold{1}}

\renewcommand{\l}{\lambda}

\newcommand{\N}{\mathbb{N}}
\renewcommand{\O}{\mathcal{O}}
\renewcommand{\P}{\mathbb{P}}
\renewcommand{\r}{\rho}
\renewcommand{\t}{\tau}
\newcommand{\R}{\mathbb{R}}
\newcommand{\s}{\sigma}

\newcommand{\Var}{\mathbb{V}\mathrm{ar}}

\DeclareMathOperator{\arcsinh}{arcsinh}
\newcommand{\barF}{\overline{F}}

\newcommand{\rhotoone}{{\r\uparrow 1}}
\newcommand{\xtoinfty}{{x\rightarrow\infty}}

\renewcommand{\bar}{\overline}
\renewcommand{\tilde}{\widetilde}

\renewcommand{\O}{O}
\newcommand{\rev}[1]{}
\newcommand{\mg}{\mathrm{M/G/1}}

\def\num{8}
\def\xstep{.04cm}
\def\xlength{40*\num-30}

\def\seqdown{(2,-2)--++(0,1)--++(1,-1)--++(0,1)--++(2,-2)--++(0,1)--++(2,-2)--++(0,1)--++(1,-1)--++(0,1) --++(2,-2)--++(0,2)--++(2,-2)--++(0,1)--++(3,-3)--++(0,1)--++(2,-2)--++(0,3)--++(2,-2)--++(0,1)--++(2,-2) --++(0,2)--++(2,-2)--++(0,1)--++(2,-2)--++(0,1)--++(1,-1)--++(0,2)--++(2,-2)--++(0,1)--++(2,-2)} 

\def\sequpshort{(0,2)--++(1,-1)--++(0,1)--++(1,-1)--++(0,2)--++(1,-1)--++(0,1)--++(1,-1)--++(0,2)--++(1,-1)
--++(0,1)--++(2,-2)--++(0,2) --++(1,-1)--++(0,2)--++(1,-1)--++(0,2)--++(1,-1)} 

\def\setupaxes{
 \draw[<->,thick] (0,100) node(yaxis)[above]{$X(t)$} |- (\xlength,0) node(xaxis)[right]{$t$};
 \node at (0,-2)[below]{$0$};
 \node at (0,85)[left]{$(1-\r)x_\r^*$};
 \node at (\xlength-25,0)[below]{$x_\r^*$};
}

\def\setuptangentbelow{
 \draw[dashed] (0,85) -- (\xlength-25,0) node(rho)[pos=0.4,sloped,below]{$\frac{dV}{dt}{ =-(1-\r)}$}; 
}

\begin{document}

\title{Uniform Asymptotics for Compound Poisson Processes with Regularly Varying Jumps and Vanishing Drift}
\author{
\begin{tabular}{ccccc}
Bart Kamphorst\textsuperscript{a} &&& Bert Zwart\textsuperscript{a,b} \\
\texttt{bart.kamphorst@cwi.nl} &&& \texttt{bert.zwart@cwi.nl}
\end{tabular} \\
\begin{tabular}{l}
\footnotesize{\textsuperscript{a} Centrum Wiskunde \& Informatica, P.O. Box 94079, 1090 GB Amsterdam, the Netherlands} \\
\footnotesize{\textsuperscript{b} Technische Universiteit Eindhoven, P.O. Box 513, 5600 MB Eindhoven, the Netherlands}
\end{tabular}
}
\date{\small{\today}}
\maketitle

\begin{abstract}
This paper addresses heavy-tailed large deviation estimates for the distribution tail of functionals of a class of spectrally one-sided L\'evy process. Our contribution is to show that these estimates remain valid in a near-critical regime. This complements recent similar results that have been obtained for the all-time supremum of such processes. Specifically, we consider local asymptotics of the all-time supremum, the supremum of the process until exiting $[0,\infty)$, the maximum jump until that time, and the time it takes until exiting $[0,\infty)$. The proofs rely, among other things, on properties of scale functions.

\end{abstract}

\textit{Keywords: compound Poisson process, $\mg$ queue, heavy traffic, large deviations, uniform asymptotics, first passage time, supremum}

\section{Introduction}
The analysis of spectrally one-sided \levy processes is a topic of fundamental interest in the stochastic processes literature \citep{kyprianou2014introductory} and arises in many applications, such as queueing \citep{debicki2012levy} and insurance risk theory \citep{asmussen1996large, embrechts1982estimates}. More generally, \levy processes and various functionals have been studied extensively over the last decades through fluctuation theory, leading to many interesting and useful 
results. If the underlying \levy measure is heavy-tailed, then exact expressions are harder to obtain and one often resorts to asymptotic estimates based on heavy-tailed large deviations. The goal of this paper is to assess the robustness of several of these approximations in a regime where the underlying \levy process has a small drift.

To make this more specific, consider the compound Poisson process with deterministic drift
\begin{equation}
 X^\r(t) := X_0 + \sum_{i=1}^{N^\r(t)} B_i - t,
\end{equation}
where $N^\r(t), t\geq 0,$ is a Poisson process\rev{REP2COM1} with a rate that depends on a drift parameter $\r$. With a slight abuse of terminology, we call $X^\r$ a compound Poisson process throughout this paper, and investigate the asymptotic behaviour of various functionals of $X^\r$ under the assumption that the i.i.d.\ nonnegative jump sizes $B_i$ have a regularly varying tail with index $\alpha>2$. The initial condition $X_0$ is equal in distribution to $B_i$ and independent of $\rho$; we present a more detailed model description in Section 2.
 The long-term drift $\E[X^\r(1)-X_0]$ of the process is negative, and of order $1-\r$. In the central limit regime, we let $\rhotoone$ so that the long-term drift tends to zero.

A functional that has received ample attention in the literature is the all-time supremum $M^\r_\infty:= \sup_{t\geq 0} X^\r(t)$. For fixed $\r$, as $x\rightarrow\infty$, the following estimate holds (see e.g.\ \citet{embrechts1982estimates, foss2007discrete, kyprianou2004ruin, maulik2006tail}):
\begin{equation}
  \P(M_\infty^\r>x) \sim \frac{\r}{\E[B_1]} (1-\r) \int_x^\infty \P(B_1>t)\dd t.
  \label{eq:htasy}
\end{equation}
This approximation can be very inaccurate when $\rhotoone$. Specifically, if $\rhotoone$ and $x=y/(1-\r)$ for fixed $y$, then $\P(M_\infty^\r>x)$ will converge to $\exp[-2(\E[B_1]/\E[B_1^2])y]$ (this is a heavy-traffic limit, cf.\ \citep{whitt1974heavy, whitt2002stochastic}).

Motivated by the contrast between these two regimes, \citet{olvera2011transition} establish an explicit threshold
\begin{equation}
  \tilde{x}_\r := \mu (\a-2) \frac{1}{1-\r}\log\frac{1}{1-\r},
  \label{eq:O-C}
\end{equation}
for some $\mu>0$ as specified in the next section, where the two regimes connect. In particular, they show that estimate \eqref{eq:htasy} remains valid when $\rhotoone$ and $x\geq (1+\e)\tilde{x}_\r$. Similar results, including examinations when the heavy-traffic approximation remains valid, can be found in \citet{blanchet2013uniform, denisov2014heavy, kugler2013upper}.

The investigation and contributions of the present paper focus on the validity of heavy-tail approximations like \eqref{eq:htasy} when $\rhotoone$.
All of the above-mentioned works focus on global asymptotics of the all-time supremum functional $M_\infty^\r$, and one may wonder how robust the obtained insights are when other functionals of importance are considered. For example, another well-studied functional of \levy processes is the first passage time of zero, $\tau^\r$, which among others may characterize a busy-period duration in queueing theory. A third functional of importance  is  $M_\tau^\r := \sup_{t<\tau^\r} X^\r(t)$. A series of prior works \citep{demeyer1980asymptotic, zwart2001tail, baltrunas2004tail} obtain useful asymptotic approximations for $\tau^\r$, while $M_\tau^\r$ has been considered in \citet{asmussen1998subexponential}. All these works focus on (a subclass of) subexponential jump sizes and fixed $\r$. Our aim is to investigate how robust these asymptotic estimates are when also $\rhotoone$.

We feel that our main achievement is a description of the tail behaviour of $\P(\tau^\r>x)$ as $\xtoinfty$ while $\rhotoone$. For fixed $\rho$, \citet{zwart2001tail} showed that 
\begin{equation}
 \P(\tau^\r>x) \sim \frac{1}{1-\r}\P(B_1>(1-\r)x)
\end{equation}
as $\xtoinfty$. In the current paper, we show that this large deviations approximation remains valid as $\rhotoone$ for all $x$ above a certain threshold $x_\r^*$ which turns out to be much larger than threshold~\eqref{eq:O-C}:
\begin{equation}
  x_\r^* := \frac{1}{(1-\r)^2} \left(\log\frac{1}{1-\r}\right)^{k^*},
  \label{eq:xrhostar}
\end{equation}
where $k^*>2$. 
We actually show that the asymptotic behaviour of $\P(\tau^\r>x)$ coincides with $\P(M_\tau^\r>(1-\r)x)$; intuitively, if the process hits zero after time $x$, then it is likely that the process obeyed the long-term drift after reaching level $(1-\r)x$ early in time. Uniform heavy-tail approximations for $M_\tau^\r$, which are also established in this paper as by-product of independent interest, yield the given asymptotic. The gap between $x_\r^*$ and $\tilde{x}_\r/(1-\r)$ is required for technical reasons; however, we show that our result does \textit{not} hold for $k=0$ in \eqref{eq:xrhostar} (i.e.\ if $x_\r^*$ is proportional to $(1-\r)^{-2}$).

Additional theorems that lead to our main result provide uniform heavy-tail approximations on the ``local'' tail probability $\P(M_\infty^\r\in[x,x+T))$ of the all-time supremum functional $M_\infty^\r$; and the largest jump $B_\tau^\r$ until time $\tau^\r$. The local asymptotics of $M_\infty^\r$ provide a generalization of Corollary 2.1(b) in \citet{olvera2011transition} and are obtained in a similar fashion via a decomposition of the so-called Pollaczek-Khinchine formula. In addition, we derive asymptotic expressions for the conditional expected time of reaching a high level $a$, given that level $a$ is reached before time $\tau^\r$. 
The corresponding lemma relies heavily on fluctuation theory for \levy\rev{REP1COM2} processes; specifically, it relies on the theory of scale functions. A recent review article on and examples of scale functions can be found in \citet{kuznetsov2013theory} and \citet{hubalek2011scale}, respectively.

The paper is organized as follows.
A precise description of the model and an introduction to the notation used can be found in Section~\ref{sec:preliminaries}. Section~\ref{sec:results} presents and discusses our results; in particular, Theorems~\ref{thm:WisNB*} and \ref{thm:PisC} display our main results. The four subsequent sections are each devoted to the proof of one theorem. Section~\ref{sec:firstpassagetime} contains the extensive proof of a crucial lemma, and, finally, Section~\ref{sec:tightness} provides the theoretical support for the discussion presented in Section~\ref{subsec:tightness}.

\section{Preliminaries}
\label{sec:preliminaries}
Let $\{B\}\cup \{B_i\}_{i=0}^\infty$ be a sequence of non-negative, independent and identically distributed (i.i.d.) regularly varying random variables \citep[cf.][]{bingham1989regular} with mean $\E[B]>0$ and finite variance $\sigma_B^2$. More specifically, their common cumulative distribution function (c.d.f.) $F_B:\R \rightarrow [0,1], F_B(0)=0$\rev{REP1COM3} is characterized by its tail
\begin{equation}
 \barF_B(x) := \P(B>x) = L(x)x^{-\a},
 \label{eq:regvarB}
\end{equation}
where $\a>2, \a\neq 3$ and $L(x)$ is a slowly varying function: $\lim_{\xtoinfty} L(ax)/L(x) = 1$ for all $a>0$. A key property of such distributions is that $\E[B^p]<\infty$ for $p<\a$ and $\E[B^p]=\infty$ for $p>\a$. The $\a$-th moment can be either finite or infinite. For technical reasons, this article does not address the $\a=3$ case. It should be noted that regularly varying distributions are a subclass of subexponential distributions \citep{goldie1998subexponential}, and as such satisfy $\lim_\xtoinfty \P(B_1+\ldots+B_n>x)/\P(B_1>x)=n$.

Define the Poisson process $N^1(t), t\geq 0$, which is independent of $B$ and has rate $1/\E[B]$\rev{REP1COM6}\rev{REP2COM3}. Then $N^\r(t):=N^1(\r t), t\geq 0$, is a Poisson process with rate $\lambda^\rho:=\rho/\E[B]$\rev{REP1COM8} and the process $X^\r:[0,\infty)\rightarrow \R$ given by
\begin{equation}
 X^\r(t) := X_0 + \sum_{i=1}^{N^\r(t)} B_i - t
\end{equation}
is a compound Poisson process with initial value $X^\r(0) = X_0 := B_0$ and long-term drift $\E[X^\r(1)-X^\r(0)]=-(1-\r)<0$. The process $X^\r(t)$ experiences a deterministic decrease of $-t$ and has jumps of size $B_i$. For this reason we refer to $F_B$ as the jump size distribution. 

The first passage time of level $x$ is denoted by $\s^\r(x):= \inf\{t\geq 0:X^\r(t)\geq x\}$, whereas the first hitting time of level zero is indicated by $\tau^\r := \inf\{t\geq 0:X^\r(t)=0\}$. 
Of primary interest in this article are the supremum $M^\r_\t$ of $X^\r(t)$ until the first down-crossing of level zero, i.e.\ $M_\tau^\r:= \sup\{X^\r(t): 0\leq t\leq \tau^\r \}$, and the all-time supremum $M_\infty^\r:= \sup\{X^\r(t): t\geq 0\}$ of the L{\'e}vy process. 
We also derive a result on the largest jump $B^\r_\tau$ before time $\tau^\r$: $B_\tau^\r:=\sup\{B_i: 0\leq i \leq N^\r(\tau^\r) \}$.

Consider the sequence of i.i.d. random variables $\{B^*\}\cup\{B_i^*\}_{i=1}^\infty$ with c.d.f.\ $F_{B^*}$. $F_{B^*}$ is the excess distribution of $B$ and will be referred to as excess jump size distribution. The excess jump size distribution can be characterized by its probability density function (p.d.f.) $f_{B^*}(x)=\frac{1}{\E[B]} \P(B>x)$ and has finite mean $\mu:=\E[B^2]/(2\E[B])<\infty$. It is assumed that $B^*$ and $B^*_i$ are independent of $N^\r,B$ and $B_i$ for all relevant indices. Since $B$ is regularly varying, Theorem 2.45 in \citet{foss2013introduction} states\rev{REP2COM4} that the tail distribution of $B^*$,
\begin{equation}
 \barF_{B^*}(x) = \frac{1}{\E[B]} \int_x^\infty \P(B>t) \dd t\sim \frac{1}{(\a-1)\E[B]}L(x) x^{-\a+1},
 \label{eq:regvarB*}
\end{equation}
is also regularly varying, where $f(z)\sim g(z)$ if and only if $\lim_{z\uparrow z^*} f(z)/g(z) = 1$ for some limiting value $z^*\in \{1,\infty\}$. In this paper, the limit of interest is either $\rhotoone, \xtoinfty$ or $a\rightarrow \infty$. The proper limit should be clear from the context. Similarly, $f(z) \gtrsim (\lesssim)\, g(z)$ denotes the relation $\liminf_{z\uparrow z^*} (\limsup_{z\uparrow z^*}) f(z)/g(z) \geq (\leq) \, 1$. We adopt the common conventions $f(z)=\O(g(z))$ if and only if $\limsup_{z\uparrow z^*} |f(z)/g(z)| < \infty$ and $f(z)=o(g(z))$ if and only if $\limsup_{z\uparrow z^*} f(z)/g(z) = 0$.\rev{REP2COM5} If both $f(z) = \O(g(z))$ and $g(z)=\O(f(z))$, then this is denoted by $f(z) = \Theta(g(z))$.

Let $T\in(0,\infty)$ be any positive constant and define the interval $\D=[0,T)$. In the remainder of this article we will denote the ``local'' probability $\P(B^*\in [x,x+T))$ by $\P(B^*\in x+\D)$. Furthermore, we adopt the well-known conventions $\lfloor x \rfloor := \max\{n\in \N: n\leq x\}$ and $\lceil x \rceil := \min\{n\in \N: n\geq x\}$.

Many expressions in this article involve constants which do not provide additional insight, and which do not contribute to the global behaviour of the expressions. For this reason, many constants have been replaced by $C$: a constant whose value may change from line to line.

Most variables that have been introduced so far depend on the parameter $\rho$. Now that their dependence has been noted, we drop the superscripts $\rho$ for the remainder of this article. Variables that are introduced in later sections and that depend on $\rho$ will have a sub- or superscript unless mentioned otherwise.

\section{Results and discussion}
\label{sec:results}

The purpose of this section is to present and discuss our main results. Our first theorem relates the local tail probability $\P(M_\infty \in x + \D)$ to the local tail probability $\P(B^*\in x + \D)$:
\begin{theorem} \label{thm:WisNB*}
Suppose $\P(B>x) = L(x)x^{-\a}$ for some $\a>2, \a\neq 3$ and $L(x)$ slowly varying. Let $\mu=\E[B^2]/(2\E[B])$ and define $x_\r:= k \mu (\a-1)\frac{1}{1-\r}\log\frac{1}{1-\r}$ for any $k>1$. Then for any fixed interval $\D=[0,T)$ the relation\rev{REP2COM6}
\begin{equation}
 \sup_{x\geq x_\r} \left|\frac{\P(M_\infty \in x+\D)}{\frac{\r}{1-\r}\P(B^*\in x+\D)} - 1\right| \rightarrow 0
 \label{eq:WisNB*}
\end{equation}
holds as $\rhotoone$. Furthermore, \eqref{eq:WisNB*} remains valid for $k=1$ provided that $L(x)/(\log x)^\a \rightarrow \infty$.
\end{theorem}

Theorem~\ref{thm:WisNB*} extends Corollary 2.3(b) of \citet{olvera2011transition}, who considered the ``global'' tail probability $\D=[0,\infty)$. The similarity of the results is also reflected in the proof of the theorem, which greatly depends on the Pollaczek-Khintchine formula and the power law nature of the jump size distribution. A key difference between the proofs is \citeauthor{olvera2011transition}'s application of the ``global'' big jump asymptotics as reported by \citet{borovkov2008asymptotic} versus our usage of the ``local'' analogues as derived by \citet{denisov2008large}. The transition point $\tilde{x}_\r$ in \citet{olvera2011transition} (cf.\ expression \eqref{eq:O-C}) differs from $x_\r$ by a factor $\frac{\a-1}{\a-2}$, which is an artefact of our analysis of the local tail probability (index $\a$) as opposed to their analysis of the global tail probability (index $\a-1$). Similarly, their $k=1$ case requires $L(x)$ to asymptotically dominate $(\log x)^{\a-1}$ instead of $(\log x)^{\a}$.

Our next result relates the tail behavior of $M_\tau$ to that of $B$:
\begin{theorem} \label{thm:CisNB}
Suppose that all conditions in Theorem~\ref{thm:WisNB*} hold. 
Then
\begin{equation}
 \sup_{x\geq x_\r} \left|\frac{\P(M_\tau>x)}{\frac{\r}{1-\r}\P(B>x)} - 1 \right| \rightarrow 0
 \label{eq:CisNB}
\end{equation}
holds as $\rhotoone$. Furthermore, \eqref{eq:CisNB} remains valid for $k=1$ provided that $L(x)/(\log x)^\a \rightarrow \infty$.
\end{theorem}
Theorem~\ref{thm:CisNB} is related to a similar result for general random walks, derived for a larger class of subexponential distributions, cf.\ \citet[Theorem~2.1]{asmussen1998subexponential}. Again, the contribution in our setting is the validity of this asymptotic estimate in the near-critical regime. Also the intuition behind this result, that $M_\tau$ is comparable in size to the largest jump $B_\tau$, remains valid:
\begin{theorem} \label{thm:BmaxisNB}
Suppose $\P(B>x) = L(x)x^{-\a}$ for some $\a>2, \a\neq 3$ and $L(x)$ slowly varying. Let $\hat{x}_\r$ satisfy $\P(B>\hat{x}_\r)/(1-\r)^2\rightarrow 0$ as $\rhotoone$. Then the relation
\begin{equation}
 \sup_{x\geq \hat{x}_\r} \left|\frac{\P(B_\tau>x)}{\frac{1}{1-\r}\P(B>x)} - 1 \right| \rightarrow 0
 \label{eq:BmaxisNB}
\end{equation}
holds as $\rhotoone$. In particular, the above statement holds for $\hat{x}_\r \geq 1/(1-\r)$.
\end{theorem}
\begin{corollary} \label{cor:CisBmax}
Suppose that all conditions in Theorem~\ref{thm:WisNB*} hold. 
Then
\begin{equation}
 \sup_{x\geq x_\r} \left|\frac{\P(M_\tau>x)}{\P(B_\tau>x)} - 1 \right| \rightarrow 0
 \label{eq:CisBmax}
\end{equation}
holds as $\rhotoone$. Furthermore, \eqref{eq:CisBmax} remains valid for $k=1$ provided that $L(x)/(\log x)^\a \rightarrow \infty$.
\end{corollary}

We are now ready to examine the asymptotic behaviour of the tail probability $\P(\tau>x)$ of the first hitting time of zero. Our result is as follows. 

\begin{theorem} \label{thm:PisC}
Suppose $\P(B>x) = L(x)x^{-\a}$ for some $\a>2,\a\neq 3$ and $L(x)$ slowly varying. For any $k^*>2$ define $x_\r^*:= \frac{1}{(1-\r)^2} \left(\log\frac{1}{1-\r}\right)^{k^*}$.\rev{REP2INTRO} 
Then both
\begin{equation}
 \sup_{x\geq x_\r^*} \bigg| \P(\tau >x \mid M_\tau>(1-\r)x) - 1 \bigg| \rightarrow 0
 \label{eq:PandCisC}
\end{equation}
and
\begin{equation}
 \sup_{x\geq x_\r^*} \bigg| \P(M_\tau>(1-\r)x \mid \tau >x) - 1 \bigg| \rightarrow 0
 \label{eq:PandCisP}
\end{equation}
hold as $\rhotoone$. In particular, \eqref{eq:PandCisC} and \eqref{eq:PandCisP} imply
\begin{equation}
 \sup_{x\geq x_\r^*} \left|\frac{\P(\tau >x)}{\P(M_\tau>(1-\r)x)} - 1\right| \rightarrow 0
 \label{eq:PandC}
\end{equation}
as $\rhotoone$.
\end{theorem}

\begin{figure}[t]
\centering
\resizebox{.9\textwidth}{!}{
\begin{tikzpicture}[x=\xstep ,y=\xstep] 
 \setupaxes
 \setuptangentbelow
 \draw (0,0) --++ \sequpshort --++ (0,\num*10+5)
 \foreach \x in {0,...,\num}
 { --++ \seqdown}
 ;
 \node at (15,90) [above]{$M_\tau$};
 \node at (35,0)[below]{$\O(1/(1-\r))$};
 \node at (280,-2)[below]{$\tau$};
\end{tikzpicture}
}
\caption{Illustration of a scenario where $X(t)$ stays positive for a long time due to a large jump early in the process. The largest jump of size $(1-\r)x$ happens at time $\O(1/(1-\r))$. The long-term drift of $-(1-\r)$ suggests that $\tau \approx (1+o(1))x$.}
\label{fig:largeCmax}
\end{figure}
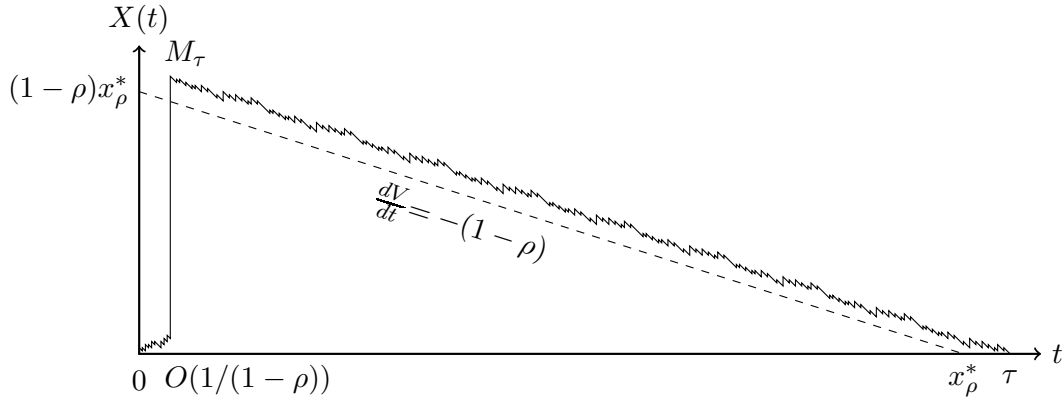

For fixed $\rho$, related questions have been examined by \citet{durrett1980conditioned} for random walks and \citet{zwart2001tail} for queues. Their results lead to the insight that a large value of $\tau$ is caused by an `early' big jump, after which the process drifts towards $0$ at rate $1-\rho$ (see Figure~\ref{fig:largeCmax}). 
%
This suggests the approximation $\tau \approx M_\tau/(1-\rho)$, which was made rigorous by \citet{zwart2001tail} using a sample-path analysis. The challenge in our setting is to show that the big jump occurs at a time that does not grow too large as $\rhotoone$. This is settled by the crucial technical Lemma~\ref{lem:LateJump} in Section~\ref{sec:hittingtime}, which essentially states that it takes $\O(1/(1-\r))$ time units until the largest jump. This lemma is proven by providing an estimate of the time until the big jump in terms of $q$-scale functions, which in turn need to be estimated in detail for various specific ranges of parameter values.

\subsection{Tightness of bounds} \label{subsec:tightness}
It is natural to question the quality of our thresholds $x_\rho$ and $x_\rho^*$ in the results presented above. The next paragraphs show that our choices are close to optimal, in the sense that our results are no longer valid if the logarithmic terms in $x_\rho$ and $x_\rho^*$ are dropped.

First, consider the function $x_\r = k\mu(\a-1)\frac{1}{1-\r} \log \frac{1}{1-\r}$ as presented in Theorem~\ref{thm:WisNB*}, Theorem~\ref{thm:CisNB} and Corollary~\ref{cor:CisBmax}. As stated earlier, Theorem~\ref{thm:WisNB*} is the local analogue of Corollary 2.3(b) in \citet{olvera2011transition} and the function $x_\r$ only differs by a constant from their function $\tilde{x}_\r$. Additionally, their Corollary 2.3(a) states that the tail probability $\P(M_\infty>x)$ asymptotically behaves as an exponential random variable for $x < (1-\e)\tilde{x}_\r$, $\e>0$ sufficiently small. This result suggests that the local tail probability $\P(M_\infty\in x+\D)$ behaves as the density of an exponential random variable for $x$ sufficiently small. The next lemma supports this suggestion by presenting a local analogue of Kingman's heavy-traffic approximation that appears to be new:

\begin{lemma} \label{lem:Wdensity}
Suppose that the jump size p.d.f.\ $f_B(x)$ of $B$ is completely monotone; i.e. $f_B(x)$ and all its derivatives exist and satisfy $(-1)^n \frac{\dd^n}{\dd x^n} f_B(x) \geq 0$ for all $x>0$ and $n=1,2,\ldots$ Fix $y>0$. Then the all-time supremum p.d.f.\ of $f_{M_\infty}(x)$ on $(0,\infty)$ exists and satisfies
\begin{equation}
 \lim_{\rhotoone} \frac{1}{1-\r}f_{M_\infty}\left(\frac{y}{1-\r}\right) = \frac{1}{\E[B^*]} e^{-\frac{y}{\E[B^*]}}.
 \label{eq:Wdensity}
\end{equation}
\end{lemma}

We hence expect $\P(M_\infty\in x+\D)$ to display exponential decay as $\rhotoone$ for $x$ sufficiently smaller than $x_\r$, similar to the results of \citet{olvera2011transition}.
 Analysing $\P(M_\infty\in x+\D)$ for general $x\leq x_\r$ is beyond the scope of this paper; however, the corollary below shows that $(1-\r)x(\r)$ must diverge to infinity in order for Theorem~\ref{thm:WisNB*} to remain true:

\begin{corollary} \label{cor:lowerboundxrho}
Suppose $\P(B>x) = L(x)x^{-\a}$ for some $\a>2$ and $L(x)$ slowly varying, and assume that the jump size p.d.f.\ $f_B(x)$ of $B$ is completely monotone. Fix $y>0$. Then for $y_\r = \frac{y}{1-\r}$ the limit
\begin{equation}
 \lim_{\rhotoone} \frac{\P(M_\infty \in y_\r + \D)}{\frac{\r}{1-\r}\P(B^* \in y_\r+ \D)}
\end{equation}
diverges to infinity.
\end{corollary}
The proof of Theorem~\ref{thm:CisNB} derives the estimates
\begin{equation}
  \frac{1}{\l}\P(M_\infty \in [x,x+1)) \lesssim \P(M_\tau > x) \lesssim \frac{1}{\l}\P(M_\infty \in [x-1,x))
  \label{eq:sandwichMtau}
\end{equation}
as $\xtoinfty$. As such, a similar necessary condition on any function $x(\r)$ for which Theorem~\ref{thm:WisNB*} holds is also necessary for  Theorem~\ref{thm:CisNB}. An analogue argument holds for Corollary~\ref{cor:CisBmax}.

We next discuss the function $x_\r^*=\frac{1}{(1-\r)^2}\left(\log \frac{1}{1-\r}\right)^k$ which is of interest in Theorem~\ref{thm:PisC}. The proof of Theorem~\ref{thm:PisC} greatly relies on Theorem~\ref{thm:CisNB} but considers $\P(M_\tau>(1-\r)x)$ instead of $\P(M_\tau>x)$. We would therefore expect Theorem~\ref{thm:PisC} to hold with $x(\r) = x_\r/(1-\r)$. The current proof, however, requires the higher level $x_\r^*$ for technical reasons. In contrast, the following lemma gives a lower bound on $x(\r)$ if it is to replace $x_\r^*$. In particular, it states that $(1-\r)^2x(\r)$ needs to diverge to infinity:

\begin{lemma} \label{lem:lowerboundxrhostar}
Suppose $\P(B>x) = L(x)x^{-\a}$ for some $\a>2$ and $L(x)$ slowly varying. Fix $y>0$. Then for $y_\r^* = \frac{y}{(1-\r)^2}$ the limit
\begin{equation}
 \lim_{\rhotoone} \frac{\P(\tau > y_\r^*)}{\frac{\r}{1-\r} \P(B > (1-\r)y_\r^*)}
\end{equation}
diverges to infinity.
\end{lemma}

\section{Local asymptotics of the all-time supremum}
\label{sec:alltimesupremum}
This section contains the proof of Theorem~\ref{thm:WisNB*}. We consider the all-time supremum by its Pollaczek-Khintchine infinite-series representation. From this representation, we distinguish between few jumps and many jumps scenarios (small and large $n$), where the threshold is approximately $x/\E[B^*]$. It is shown that under the few jumps scenario, a large all-time supremum is most probably due to a large value of a single $B_i^*$. Contrastingly, the many jumps scenario is shown to be negligible.

Define $S^*_0:=0$ and $S^*_n:= \sum_{i=1}^n B_i^*$. By Theorem VIII.5.7 in \citet{asmussen2003applied},
\begin{equation}
\P(M_\infty \in x+\D) = \sum_{n=0}^\infty (1-\r)\r^n\P(S^*_n\in x + \D) \label{eq:WisS}
\end{equation}
for all $x>0$. An equivalent representation of \eqref{eq:WisNB*} is therefore
\begin{equation}
 \sup_{x\geq x_\r} \left|\frac{\sum_{n=1}^\infty (1-\r)\r^n \left[\P(S^*_n \in x+\D)-n\P(B^*\in x+\D)\right]}{\frac{\r}{1-\r}\P(B^*\in x+\D)}\right| \rightarrow 0
 \label{eq:SisnBx}
\end{equation}
as $\rhotoone$. Fix a constant $\d$ that satisfies $\max\{\frac{1}{2},\frac{1}{\a-1}\}<\d<1$ and define $U_\d(x):=\lfloor (x-x^\d)/\mu\rfloor$. Then the numerator in~\eqref{eq:SisnBx} can be decomposed as
\begin{align*}
 \left|\sum_{n=1}^\infty (1-\r)\r^n \left[\P(S^*_n \in x+\D)-n\P(B^*\in x+\D)\right]\right| \hspace{-95pt} & \\
  &\leq \sum_{n=1}^{U_\d(x)} (1-\r)\r^n\bigg|\P(S^*_n\in x + \D) - n\P(B^*\in x-(n-1)\mu+\D)\bigg| \\
   &\quad + \sum_{n=1}^{U_\d(x)} (1-\r)\r^n n \bigg|\P(B^*\in x -(n-1)\mu + \D) - \P(B^*\in x+\D)\bigg| \\
    &\quad + \left|\sum_{n=U_\d(x)+1}^\infty (1-\r)\r^n\left[\P(S^*_n\in x + \D) - n\P(B^*\in x+\D)\right]\right|. \numberthis \label{eq:WisSdecomp}
\end{align*}
Here, the first term corresponds to the few jumps scenario and the third term corresponds to the many jumps scenario. The second term corrects a shift in the argument of $\P(B^*\in \cdot)$, which is required for application of the following lemma: 

\begin{lemma} \label{lemma:DDS}
 Suppose $\xi$ is a non-negative regularly varying random variable whose c.d.f.\ has index $-\a_\xi<-2, \a_\xi \neq 3$; i.e.\ $\P(\xi>x)=L(x)x^{-\a_\xi}$. Let $F_{\xi^*}$ be the excess distribution of $\xi$ with index $-\a_\xi+1<-1$ and i.i.d. samples $\xi^*,\xi^*_1,\xi^*_2,\ldots$ For any $\max\left\{\frac{1}{\a_\xi-1},\frac{1}{2}\right\}<\Gamma<1$ denote $U_\Gamma(x)=\lfloor\frac{x-x^\Gamma}{\E[\xi^*]}\rfloor$. Then, there exists a non-increasing function $\phi(x)$ satisfying $\phi(x)\downarrow 0$ as $\xtoinfty$ such that
 \begin{equation*}
  \sup_{1\leq n\leq U_\Gamma(x)\rev{REP2COM10}} \left|\frac{\P(\xi^*_1+\ldots+\xi^*_n \in x+\D)}{n\P(\xi^*\in x-(n-1)\E[\xi^*]+\D)}-1\right| \leq \phi(x). 
 \end{equation*}
\end{lemma}
The proof is delayed until the end of this section and relies heavily on the machinery provided by \citet{denisov2008large}.
Lemma~\ref{lemma:DDS} is closely related to the property $\lim_{\xtoinfty} \P(B_1^*+\ldots+B_n^*>x)/\P(B_1^*>x)=n$ and guarantees that, for some non-increasing $\phi(x)\downarrow 0$ as $\xtoinfty$, expression \eqref{eq:WisSdecomp} is dominated by
\begin{align*}
  &\phi(x) \sum_{n=1}^{U_\d(x)} (1-\r)\r^nn\P(B^*\in x+\D) \\
    &\quad + (1+\phi(x))\sum_{n=1}^{U_\d(x)} (1-\r)\r^n n\bigg|\P(B^*\in x -(n-1)\mu + \D) - \P(B^*\in x+\D)\bigg| \\
    &\qquad + \sum_{n=U_\d(x)+1}^\infty (1-\r)\r^n\left[1 + n\P(B^*\in x+\D)\right] \\
  &= \phi(x) \mathrm{I} + (1+\phi(x))\mathrm{II} + \mathrm{III}
\end{align*}
Term $\mathrm{I}$ is bounded by $\frac{\r}{1-\r}\P(B^*\in x+\D)$, so that $x_\r\rightarrow \infty$ implies
\begin{equation}
 \sup_{x\geq x_\r} \frac{\phi(x) \mathrm{I}}{\frac{\r}{1-\r}\P(B^*\in x+\D)} \leq \phi(x_\r) \rightarrow 0 \label{eq:SisnBx:I}
\end{equation}
as $\rhotoone$.

Error term $\mathrm{II}$ is split into two parts. Fix $\g$ such that $0<\g<\d$ 
and define $V_\g(x):= \lfloor (1-\g)x/\mu\rfloor$. For $x$ sufficiently large we have $V_\g(x)<U_\d(x)$ and $\mathrm{II}$ may be written as
\begin{align*}
 \sum_{n=1}^{U_\d(x)} & (1-\r)\r^n n \bigg|\P(B^*\in x -(n-1)\mu + \D) - \P(B^*\in x+\D)\bigg| \\
  &= \sum_{n=1}^{V_\g(x)} (1-\r)\r^n n \bigg|\P(B^*\in x -(n-1)\mu + \D) - \P(B^*\in x+\D)\bigg| \\
   &\qquad + \sum_{n=V_\g(x)+1}^{U_\d(x)} (1-\r)\r^n n \bigg|\P(B^*\in x -(n-1)\mu + \D) - \P(B^*\in x+\D)\bigg| \\
  &= \mathrm{IIa} + \mathrm{IIb}.
\end{align*}
For $1\leq n\leq V_\g(x)$, Newton's generalized binomial Theorem implies that
\begin{align*}
 \frac{\P(B^*\in x-(n-1)\mu +\D)}{\P(B^*\in x+\D)} \hspace{-12pt}&\hspace{12pt} \leq \frac{\P(B>x-(n-1)\mu)}{\P(B>x+T)}
  \sim \left(1-\frac{(n-1)\mu+T}{x+T}\right)^{-\a} \\
  &= 1+\sum_{m=1}^\infty \frac{\a(\a+1)\cdots(\a+m-1)}{m!} \left(\frac{(n-1)\mu+T}{x+T}\right)^m \\
  &\leq 1+\a\left(\frac{(n-1)\mu+T}{x+T}\right)\left(1-\frac{(n-1)\mu+T}{x+T}\right)^{-\a-1} \\
  &\lesssim 1+\a\g^{-\a-1}\frac{(n-1)\mu+T}{x+T},
\end{align*}
as $\xtoinfty$.\rev{REP2COM11} Therefore,
\begin{equation*}
 \frac{\P(B^*\in x-(n-1)\mu +\D)}{\P(B^*\in x+\D)} - 1 \lesssim C \frac{n-1}{x}
\end{equation*}
as $\xtoinfty$. Substituting this into $\mathrm{IIa}$ gives
\begin{align*}
  \mathrm{IIa} &\lesssim C\P(B^*\in x+\D) \frac{1}{x}\sum_{n=1}^{V_\g(x)} (1-\r)\r^n n (n-1) \\ 
    &\leq C\P(B^*\in x+\D) \frac{2\r^2}{(1-\r)^2x}\left(1-\r^{V_\g(x)}-(1-\r)V_\g(x)\r^{V_\g(x)}\right)
    \leq \frac{C\r}{(1-\r)^2x}\P(B^*\in x+\D),
\end{align*}
and hence
\begin{equation}
 \sup_{x\geq x_\r} \frac{\mathrm{IIa}}{\frac{\r}{1-\r}\P(B^*\in x+\D)} \lesssim \frac{C}{\log\frac{1}{1-\r}}\rightarrow 0 \label{eq:SisnBx:IIa}
\end{equation}
as $\rhotoone$.

Next, we consider term $\mathrm{IIb}$. Since $\P(B^*\in y+\D)$ is decreasing in $y$, term $\mathrm{IIb}$ is bounded by
\begin{align*}
 \mathrm{IIb} &\leq C(1-\r)\r^{V_\g(x)+1} x \sum_{n=V_\g(x)+1}^{U_\d(x)} \P(B^*\in x-(n-1)\mu +\D) \\
  &\leq C(1-\r)\r^{(1-\g)x/\mu} x \int_{x-\mu U_\d(x)}^{x-\mu V_\g(x)} \P(B^*\in t+\D) \dd t. 
\end{align*}
Noting that $\P(B^*\in x+\D)$ is regularly varying with index $-\a<-2$, Theorem 1.5.11 in \citep{bingham1989regular} indicates that
\begin{align*}
 \mathrm{IIb} &\lesssim C(1-\r)\r^{(1-\g)x/\mu} x (x-\mu U_\d(x)) \P(B^*\in x-\mu U_\d(x)+\D) \\
  &\leq C(1-\r)\r^{(1-\g)x/\mu} x^{1+\d} \P(B^*\in x^\d-\D).
\end{align*}
It remains to verify that $\mathrm{IIb}$ decreases sufficiently fast for $x\geq x_\r$. One can see that
\begin{align*}
 \sup_{x\geq x_\r} \frac{\mathrm{IIb}}{\frac{\r}{1-\r}\P(B^*\in x+\D)} &\lesssim C \sup_{x\geq x_\r} (1-\r)^2\r^{(1-\g)\frac{x}{\mu}-1} x^{1+\d} \frac{\P(B^*\in x^\d+\D)}{\P(B^*\in x+\D)} \\
  &\leq C \sup_{x\geq x_\r} (1-\r)^2\r^{(1-\g)\frac{x}{\mu}-1} x^{1+\d} \frac{\P(B>x^\d)}{\P(B>x+T)} \\
  &\sim C \sup_{x\geq x_\r} (1-\r)^2 e^{\left((1-\g)\frac{x}{\mu}-1\right)\log \r} x^{1+\d+(1-\d)\a} \\
  &\leq C \sup_{x\geq x_\r} (1-\r)^2 e^{-\left((1-\g)\frac{x}{\mu}-1\right)(1-\r)} x^{1+\d+(1-\d)\a}, 
\end{align*}
where we have used $\log \r \leq -(1-\r)$ for all $\r\geq 0$. 
%
Additionally, for $\r$ sufficiently close to one, the supremum is achieved in $x=x_\r$ and
\begin{equation*}
 \sup_{x\geq x_\r} \frac{\mathrm{IIb}}{\frac{\r}{1-\r}\P(B^*\in x+\D)} \lesssim C(1-\r)^2 e^{-(1-\g)(1-\r)\frac{x_\r}{\mu}} x_\r^{1+\d+(1-\d)\a}.
\end{equation*}
Substituting $x_\r=k \mu(\a-1)\frac{1}{1-\r}\log\frac{1}{1-\r}$ now gives
\begin{align*}
 \sup_{x\geq x_\r} \frac{\mathrm{IIb}}{\frac{\r}{1-\r}\P(B^*\in x+\D)} &\lesssim C(1-\r)^{k(1-\g)(\a-1)-(1-\d)(\a-1)} \left(\log \frac{1}{1-\r}\right)^{1+\d+(1-\d)\a}
  \rightarrow 0 \numberthis \label{eq:SisnBx:IIb}
\end{align*}
as $\rhotoone$, since $\g<\d$
. This verifies the convergence of term $\mathrm{II}$.

We continue with the analysis of term $\mathrm{III}$. This term is\rev{REP2COM12} rewritten into two smaller terms:
\begin{align*}
 \mathrm{III} &= \r^{U_\d(x)+1} + \left[\left(U_\d(x)+1\right)\r^{U_\d(x)+1}+\frac{\r^{U_\d(x)+2}}{1-\r}\right]\P(B^*\in x+\D) \\
  &\leq \r^{\frac{x-x^\d}{\mu}} + \left[\left(\frac{x-x^\d}{\mu}+1\right)+ \frac{\r}{1-\r}\right]\P(B^*\in x+\D)\r^{\frac{x-x^\d}{\mu}} \\
  &\leq \r^{\frac{x-x^\d}{\mu}} + C \frac{(1-\r)x+1}{1-\r}\P(B^*\in x+\D)\r^{\frac{x-x^\d}{\mu}} \\
  &= \mathrm{IIIa} + \mathrm{IIIb}.
\end{align*}
We consider terms $\mathrm{IIIa}$ and $\mathrm{IIIb}$ in order.

For term $\mathrm{IIIa}$, we first assume that $k>1$. Potter's Theorem (e.g.\ Theorem 1.5.6 in \citet{bingham1989regular}) suggests that $\P(B^*\in x+\D)\geq T\P(B\geq x+T) \geq T C (x+T)^{-\a-\nu}$ for any fixed $\nu>0$ and $x$ sufficiently large. 
In particular, for $0<\nu < (k-1)(\a-1)$,
\begin{equation*}
 \sup_{x\geq x_\r} \frac{\mathrm{IIIa}}{\frac{\r}{1-\r}\P(B^*\in x+\D)} \leq \sup_{x\geq x_\r} C (1-\r)(x+T)^{\a+\nu}\r^{\frac{x-x^\d}{\mu}-1}.
\end{equation*}
%
Again, the supremum is achieved in $x=x_\r$ for $\r$ sufficiently close to one and hence
\begin{align*}
 \sup_{x\geq x_\r} \frac{\mathrm{IIIa}}{\frac{\r}{1-\r}\P(B^*\in x+\D)} &\leq C(1-\r)e^{(\a+\nu)\log x_\r+\left(\frac{x_\r-x_\r^\d}{\mu}-1\right)\log\r+(\a+\nu)\log\left(1+\frac{T}{x_\r}\right)} \\
   &\leq Ce^{(\a+\nu)\log x_\r-\left(\frac{x_\r-x_\r^\d}{\mu}-1\right)(1-\r)-\log\frac{1}{1-\r}+(\a+\nu)\log\left(1+\frac{T}{x_\r}\right)}.
\end{align*}
Substitution of $x_\r=k\mu(\a-1)\frac{1}{1-\r}\log\frac{1}{1-\r}$ now yields
\begin{equation}
 \sup_{x\geq x_\r} \frac{\mathrm{IIIa}}{\frac{\r}{1-\r}\P(B^*\in x+\D)} \leq C(1-\r)e^{(\a+\nu-1)\log \frac{1}{1-\r}-k(\a-1)\log\frac{1}{1-\r}+o\left(\log\frac{1}{1-\r}\right)},
 \label{eq:SisnBx:IIIa}
\end{equation}
which tends to zero as $\rhotoone$.

Alternatively, assume $k=1$ and $L(x)/(\log x)^\a\rightarrow \infty$. Then, there exists a non-increasing function $\phi(x)\downarrow 0$ such that $L(x) \geq (\log^\a x)/\phi(x)$. Similar to the preceding analysis we find
\begin{align*}
 \sup_{x\geq x_\r} \frac{\mathrm{IIIa}}{\frac{\r}{1-\r}\P(B^*\in x+\D)} &\leq \sup_{x\geq x_\r} \frac{1}{L(x+T)} (1-\r)(x+T)^{\a}\r^{\frac{x-x^\d}{\mu}-1} \\
  &\lesssim \phi(x_\r) \sup_{x\geq x_\r} \frac{1}{\log^\a x} e^{\a\log x + \left(\frac{x-x^\d}{\mu}-1\right)\log \rho - \log \frac{1}{1-\r}} \\
  &\leq \phi(x_\r) \frac{1}{\log^\a x_\r} e^{(\a-1)\log \frac{1}{1-\r} + \a\log\log \frac{1}{1-\r} - (\a-1)\left(1-x_\r^{\d-1}-x_\r^{-1}\right)\log \frac{1}{1-\r}} \\
  &= C \phi(x_\r) \frac{1}{\log^\a \frac{1}{1-\r}} e^{\a\log \log \frac{1}{1-\r} + (\a-1)\left(x_\r^{\d-1}+x_\r^{-1}\right)\log \frac{1}{1-\r}} \\
  &= C \phi(x_\r) e^{(\a-1)\left(x_\r^{\d-1}+x_\r^{-1}\right)\log \frac{1}{1-\r}} \rightarrow 0
\end{align*}
as $\rhotoone$ since $(\log x)/x^{1-\d}\rightarrow 0$ for any $\d<1$.

Finally, for term $\mathrm{IIIb}$ one can see that
\begin{equation*}
 \sup_{x\geq x_\r} \frac{\mathrm{IIIb}}{\frac{\r}{1-\r}\P(B^*\in x+\D)} = C\sup_{x\geq x_\r} ((1-\r)x+1)\r^{\frac{x-x^\d}{\mu}-1}.
\end{equation*}
As before, the supremum is attained in $x=x_\r$ for $\r$ sufficiently close to one. Thus,
\begin{align*}
 \sup_{x\geq x_\r} \frac{\mathrm{IIIb}}{\frac{\r}{1-\r}\P(B^*\in x+\D)} &= C((1-\r)x_\r+1) e^{\left(\frac{x_\r-x_\r^\d}{\mu}-1\right) \log \rho} \\
  &\leq C\log\frac{1}{1-\r} e^{-k(\a-1)(1+o(1))\log\frac{1}{1-\r}}
  \rightarrow 0 \numberthis \label{eq:SisnBx:IIIb}
\end{align*}
as $\rhotoone$.
From \eqref{eq:SisnBx:I}, \eqref{eq:SisnBx:IIa}, \eqref{eq:SisnBx:IIb}, \eqref{eq:SisnBx:IIIa} and \eqref{eq:SisnBx:IIIb}, we may conclude that \eqref{eq:SisnBx} and equivalently \eqref{eq:WisNB*} converges to zero. This completes the proof of Theorem~\ref{thm:WisNB*}.

The section is concluded by the proof of Lemma~\ref{lemma:DDS}.

\subsection{Proof of Lemma~\ref{lemma:DDS}}
First consider the case $-\a_\xi<-3$. Then $\sigma_{\xi^*}^2=\Var(\xi^*)=\frac{\E[\xi^3]}{3\E[\xi]}$ is finite, and therefore $\bar{\xi}_i^* = \frac{\xi^*_i-\E[\xi^*]}{\sigma_{\xi^*}}$ and $\bar{S}_n^* = \frac{\xi^*_1+\ldots+\xi^*_n-n\E[\xi^*]}{\sigma_{\xi^*}}$ are well-defined for all $i\geq 1, n\geq 1$. Since
\begin{equation}
 \frac{\P(\xi^*_1 + \ldots + \xi^*_n \in x+\D)}{n\P(\xi^*\in x-(n-1)\E[\xi^*]+\D)} = \frac{\P\left(\bar{S}_n^* \in \frac{x-n\E[\xi^*]+\D}{\sigma_{\xi^*}}\right)}{n\P\left(\bar{\xi}^*_1 \in \frac{x-n\E[\xi^*]+\D}{\sigma_{\xi^*}}\right)},
\end{equation}
the result follows from Theorem 8.1 in \citet{denisov2008large} once we show that \\
$(x-n\E[\xi^*])/\sqrt{(\a_\xi-3)n\log n}\rightarrow \infty$ uniformly for $1\leq n\leq U_\Gamma(x)$ as $\xtoinfty$. As $\Gamma>\frac{1}{2}$, it follows that
\begin{equation*}
 \frac{x-n\E[\xi^*]}{\sqrt{(\a_\xi-3)n\log n}}\geq \frac{x-U_\Gamma(x)\E[\xi^*]}{\sqrt{(\a_\xi-3)U_\Gamma(x)\log U_\Gamma(x)}} \sim \sqrt{\frac{\E[\xi^*]}{\a_\xi-3}}\frac{x^{\Gamma-\frac{1}{2}}}{\log x}
\end{equation*}
indeed tends to infinity as $\xtoinfty$.\rev{REP1COM9}

Now assume $-3<-\a_\xi<-2$. Let $\tilde{\xi}_i^* = \xi^*_i-\E[\xi^*]$ and $\tilde{S}_n^* = \xi^*_1+\ldots+\xi^*_n-n\E[\xi^*]$ for all $i\geq 1, n\geq 1$. Then
\begin{equation}
 \frac{\P(\xi^*_1 + \ldots + \xi^*_n \in x+\D)}{n\P(\xi^*\in x-(n-1)\E[\xi^*]+\D)} = \frac{\P\left(\tilde{S}_n^* \in x-n\E[\xi^*]+\D\right)}{n\P\left(\tilde{\xi}^*_1 \in x-n\E[\xi^*]+\D\right)}.
\end{equation}
Fix $\Gamma^*$ such that $\frac{1}{\a_\xi-1}<\Gamma^*<\Gamma$. Theorem 9.1 in \citet{denisov2008large} implies that $\P(\tilde{S}^*_n \in x+\D)\sim n\P(\tilde{\xi}^*_1\in x+\D)$ uniformly for $x\geq n^{\Gamma^*}$. The proof is concluded by showing that $(x-n\E[\xi^*])/n^{\Gamma^*} \rightarrow \infty$ uniformly for $1\leq n\leq U_\Gamma(x)$, which follows from
\begin{equation*}
 \frac{x-n\E[\xi^*]}{n^{\Gamma^*}}\geq \frac{x-U_\Gamma(x)\E[\xi^*]}{U_\Gamma(x)^{\Gamma^*}} \sim \E[\xi^*]^{\Gamma^*} x^{\Gamma-\Gamma^*}.
\end{equation*}


\section[Asymptotics of the supremum Mtau]{Asymptotics of the supremum $M_\tau$}
\label{sec:stoppingtimesupremum}
This section is dedicated to the proof of Theorem~\ref{thm:CisNB}. \citet[Section 29]{takacs1967combinatorial} 
and \citet{cohen1968extreme} 
have independently shown that
\begin{equation*}
 \P(M_\tau>x) = \frac{1}{\l} \frac{d}{dy} \log \P(M_\infty < y)\big|_{y=x},
 \label{eq:Mtau}
\end{equation*}
of which we analyse the right-hand side by means of Theorem~\ref{thm:WisNB*} and the theory of scale functions. The definition and some properties of scale functions are provided in Section~\ref{sec:firstpassagetime}; for now we only state that $\P(M_\infty < x) = (1-\r)W_\r(x)$ for a scale function $W_\r(x)$ satisfying $W_\r(x)>0$ for $x>0$. A property of interest for the current section is that $\frac{d}{dy}\log W_\r(y)$ is non-increasing and positive (cf.\ equation \ref{eq:logW}), and hence $\log W_\r(x)$ is concave.

Rewriting~\eqref{eq:Mtau} in terms of the scale function $W_\r(x)$ and exploiting its concavity gives
\begin{align*}
 \P(M_\tau>x) 
  &= \frac{1}{\l}\frac{d}{dy}\log W_\r(y)\big|_{y=x}
  \leq \frac{1}{\l}\left[\log W_\r(x)-\log W_\r(x-1)\right] \numberthis \label{eq:CisNB:bound}  \\
  &= \frac{1}{\l}\log \frac{W_\r(x)}{W_\r(x-1)}
  \leq \frac{1}{\l}\left[\frac{W_\r(x)}{W_\r(x-1)}-1\right]
  = \frac{1}{\l}\frac{\P(M_\infty\in [x-1,x))}{\P(M_\infty < x-1)}
\end{align*}
for all $x\geq 1$\rev{REP1COM12}. Theorem~\ref{thm:WisNB*} then implies
\begin{equation*}
 \P(M_\tau>x) \lesssim \frac{1}{\l}\frac{1}{\P(M_\infty < x-1)}\frac{\r}{1-\r}\P(B^*\in[x-1,x))
\end{equation*}
for $x\geq x_\r$. Applying the simple bound $\P(B^*\in[x-1,x)) \leq \frac{1}{\E[B]}\P(B > x-1)$ yields
\begin{equation*}
 \P(M_\tau>x) \lesssim \frac{1}{\P(M_\infty < x-1)} \frac{1}{1-\r}\P(B > x-1),
\end{equation*}
which, since $B$ is long-tailed, is asymptotically equivalent to
\begin{equation}
 \P(M_\tau>x) \lesssim \frac{1}{\P(M_\infty < x-1)} \frac{1}{1-\r}\P(B > x). \label{eq:CisNB:upperbound}
\end{equation}
This concludes the upper bound analysis for $\P(M_\tau>x)$.
Similarly, $\P(M_\tau>x)$ can be bounded from below. Using the inequality $\log x\geq 1-\frac{1}{x}$ for all $x\geq 0$ and slightly altering the above analysis dictates
\begin{equation}
 \P(M_\tau>x) \gtrsim \frac{1}{\P(M_\infty < x+1)} \frac{1}{1-\r}\P(B > x) \geq \frac{\r}{1-\r}\P(B > x)\rev{REP2COM14} \label{eq:CisNB:lowerbound}
\end{equation}
for all $x\geq x_\r$ as $\rhotoone$. Combining both bounds completes the proof.

\section{Asymptotics of the supremum jump size}
\label{sec:supremumjump}
This section contributes the proof of Theorem~\ref{thm:BmaxisNB}. The following equality is an interpretation of expression (3.4) in \citet{boxma1978longest}:
\begin{equation}
\P(B_\tau>x) = \P(B>x) + \int_0^x\left[1-e^{-\l\P(B_\tau>x)t}\right]\dd \P(B\leq t).\rev{REP1COM14}\rev{REP2COM15}
\end{equation}
From this equality it follows that
\begin{align*}
 \frac{\P(B>x)}{\P(B_\tau>x)} &= 1 - \int_0^x\left[\frac{1-e^{-\l\P(B_\tau>x)t}}{\P(B_\tau>x)}\right]\dd \P(B\leq t) \\
  &= 1 - \l\int_0^x t\dd \P(B\leq t) + \l \int_0^x\left[1-\frac{1-e^{-\l\P(B_\tau>x)t}}{\l\P(B_\tau>x)t}\right]t \dd \P(B\leq t) \\
  &= 1 - \r + \l\int_x^\infty t\dd \P(B\leq t) + \l \int_0^x\left[1-\frac{1-e^{-\l\P(B_\tau>x)t}}{\l\P(B_\tau>x)t}\right]t \dd \P(B\leq t),
\end{align*}
so that
\begin{multline}
 \frac{1}{1-\r}\frac{\P(B>x)}{\P(B_\tau>x)} - 1 = \frac{\l}{1-\r}\int_x^\infty t\dd \P(B\leq t) \\
  + \frac{\l}{1-\r} \int_0^x\left[1-\frac{1-e^{-\l\P(B_\tau>x)t}}{\l\P(B_\tau>x)t}\right]t \dd \P(B\leq t).
  \label{eq:BdivBmax}
\end{multline}
Note that the right-hand side of the latter expression is non-negative because $(1-e^{-y})/y\leq 1$.

The first integral in \eqref{eq:BdivBmax} can be upper bounded as
\begin{align*}
  \sup_{x\geq \hat{x}_\r} \frac{\l}{1-\r}\int_x^\infty t\dd \P(B\leq t) &= \sup_{x\geq \hat{x}_\r} \frac{\l}{1-\r}\E\left[B\,\ind\{B>x\}\right] \\
    &= \sup_{x\geq \hat{x}_\r} \left( \frac{\l}{1-\r}\E\left[B-x\bigg| B>x\right]\P(B>x) + \frac{\l x\P(B>x)}{1-\r} \right) \\
    &\lesssim C \sup_{x\geq \hat{x}_\r} \frac{\l x\P(B>x)}{1-\r}, \numberthis \label{eq:BdivBmax1}
\end{align*}
since $\E\left[B-x\mid B>x\right]\sim \frac{x}{\a-1}$, as shown in \citet[p.162]{embrechts1997modelling}\rev{REP1COM15}. Clearly, this upper bound tends to zero for all $x\geq \hat{x}_\r$ provided that $\sup_{x\geq \hat{x}_\r} \frac{x \P(B>x)}{1-\r}\rightarrow 0$ as $\rhotoone$.

Sequentially, we consider the second integral in \eqref{eq:BdivBmax}. The bound $e^y\geq 1+y+y^2/2$ for $y\geq 0$ implies
\begin{align*}
 \sup_{x\geq \hat{x}_\r} \frac{\l}{1-\r} \int_0^x\left[1-\frac{1-e^{-\l\P(B_\tau>x)t}}{\l\P(B_\tau>x)t}\right]t \dd \P(B\leq t) \hspace{-105pt} & \\
   &\leq \sup_{x\geq \hat{x}_\r} \frac{\l}{1-\r} \int_0^x\left[1-\frac{1-\frac{1}{1 + \l\P(B_\tau>x)t + \l^2\P(B_\tau>x)^2t^2/2}}{\l\P(B_\tau>x)t}\right]t \dd \P(B\leq t) \\
   &= \sup_{x\geq \hat{x}_\r} \frac{\l}{2}\frac{1}{1-\r} \int_0^x\frac{\l\P(B_\tau>x)t + \l^2\P(B_\tau>x)^2t^2}{1 + \l\P(B_\tau>x)t + \l^2\P(B_\tau>x)^2t^2/2} t \dd \P(B\leq t) \\
   &\leq \sup_{x\geq \hat{x}_\r} \frac{\l^2}{2}\frac{\P(B_\tau>x)}{1-\r} \int_0^x \left[t^2 + \l\P(B_\tau>x)t^3\right] \dd \P(B\leq t) \\
   &\leq \sup_{x\geq \hat{x}_\r} \frac{\l^2}{2}\frac{\P(B_\tau>x)}{1-\r} \left[1 + \l\P(B_\tau>x)x\right] \E[B^2].
\end{align*}
Equation \eqref{eq:BdivBmax} implies $\P(B_\tau>x)\leq \frac{\P(B>x)}{1-\r}$, so that
\begin{align*}
  \sup_{x\geq \hat{x}_\r} \frac{\l}{1-\r} \int_0^x\left[1-\frac{1-e^{-\l\P(B_\tau>x)t}}{\l\P(B_\tau>x)t}\right]t \dd \P(B\leq t) 
    &\leq \sup_{x\geq \hat{x}_\r} \frac{\l^2}{2}\frac{\P(B>x)}{(1-\r)^2} \left[1 + \frac{\l x\P(B>x)}{1-\r}\right] \E[B^2].
  \numberthis \label{eq:BdivBmax2}
\end{align*}
The second expression hence vanishes for all $x\geq \hat{x}_\r$ if both $\sup_{x\geq \hat{x}_\r} \frac{x \P(B>x)}{1-\r}$ and $\sup_{x\geq \hat{x}_\r} \frac{\P(B>x)}{(1-\r)^2}$ tend to zero as $\rhotoone$.
 From \eqref{eq:BdivBmax}, \eqref{eq:BdivBmax1} and \eqref{eq:BdivBmax2}, it follows that
\begin{equation}
 \sup_{x\geq \hat{x}_\r} \left|\frac{\P(B_\tau>x)}{\frac{1}{1-\r}\P(B>x)} - 1\right| \rightarrow 0
\end{equation}
as $\rhotoone$, provided $\sup_{x\geq \hat{x}_\r} \frac{x \P(B>x)}{1-\r}\rightarrow 0$ and $\sup_{x\geq \hat{x}_\r} \frac{\P(B>x)}{(1-\r)^2}\rightarrow 0$. These conditions are analysed by Potter's Theorem, which states that for any $0<\nu<\a-2$ there exists a constant $C_\nu>0$ such that $\P(B>x) \leq C_\nu x^{-\a+\nu}$ for all $x$ sufficiently large. In particular, we find
\begin{equation*}
  \sup_{x\geq 1/(1-\r)} \frac{x \P(B>x)}{1-\r} \leq C_\nu \sup_{x\geq 1/(1-\r)} \frac{x^{1-\a+\nu}}{1-\r} \rightarrow 0,
\end{equation*}
and similarly $\sup_{x\geq 1/(1-\r)} \frac{\P(B>x)}{(1-\r)^2} \rightarrow 0$, implying that the theorem holds for $\hat{x} \geq 1/(1-\r)$. Since $\frac{x\P(B>x)}{1-\r} \leq \frac{\P(B>x)}{(1-\r)^2}$ when $x\leq 1/(1-\r)$, the condition $\sup_{x\geq \hat{x}_\r} \frac{\P(B>x)}{(1-\r)^2}\rightarrow 0$ implies \eqref{eq:BmaxisNB} for all $\hat{x}_\r \leq \frac{1}{1-\r}$.

\section{Asymptotics of the first hitting time of level zero}
\label{sec:hittingtime}
This section is devoted to the proof of\rev{REP1COM24} Theorem~\ref{thm:PisC}. We will validate expression~\eqref{eq:PandCisC}, which considers the asymptotic behaviour of the conditional probability $\P(\tau>x \mid M_\tau>(1-\r)x)$, and expression~\eqref{eq:PandC}, which considers the asymptotic behaviour of the unconditional probability $\P(\tau>x)$ as $\rhotoone$. Expressions \eqref{eq:PandCisC} and \eqref{eq:PandC} together imply expression \eqref{eq:PandCisP} through the inequality 
\begin{equation}
  \bigg|\P(Q\mid R)-1 \bigg|
    \leq \bigg|\P(R\mid Q)-1 \bigg|\times \left|\frac{\P(Q)}{\P(R)}\right|
      + \left|\frac{\P(Q)}{\P(R)}-1\right|.
\end{equation}
for two events $Q$ and $R$ of non-zero probability.

Section~\ref{subsec:PisClowerbound} validates the asymptotic behaviour of $\P(\tau>x \mid M_\tau>(1-\r)x)$. Thereafter, Section~\ref{subsec:PisCupperbound} proves the asymptotic behaviour of $\P(\tau>x)$ by means of a sample-path analysis that makes a distinction based on the supremum $M_\tau$. The resulting events are again discriminated based on the number of jumps before $\tau$ or the first passage time of a specific level.

\subsection{Asymptotics of conditional first hitting time}
\label{subsec:PisClowerbound}
We first prove expression~\eqref{eq:PandCisC}. Since $\P(\tau>x \mid M_\tau>(1-\r)x)\geq 1$, we only need to show that $\sup_{x\geq x_\r^*} \P(\tau>x \mid M_\tau>(1-\r)x) - 1 \geq 0$ as $\rhotoone$.

Fix $p\in\left(\frac{1}{2}+\frac{1}{k^*},1\right)$ and define
$h_u(x,\r):= (1-\r)x + g(x,\r)$,
where
$g(x,\rho) := (1-\r)^{2p-1}x^p$.
The function $h_u(x,\r)$ is an upper bound for, yet asymptotically equivalent to, $(1-\r)x$. By conditioning on the event $\{M_\tau > h_u(x,\r)\}$, the long term drift $-(1-\r)$ of $X(t)$ implies that $\P(\tau>x \mid M_\tau > h_u(x,\r))$ must tend to one.

To make this precise we follow the proof of Proposition 3.1 of \citet{zwart2001tail}. Noting that $\{\sigma(y)<\tau\} = \{M_\tau>y\}$, the joint probability $\P(\tau>x; M_\tau>(1-\r)x)$ is lower bounded as
\begin{align*}
  \P(\tau>x; M_\tau>(1-\r)x) &\geq \P(\tau>x; M_\tau>h_u(x,\r)) \\
    &\geq \P(\tau - \sigma(h_u(x,\r)) > x \mid  \sigma(h_u(x,\r)) < \tau) \P(M_\tau> h_u(x,\r)),
\end{align*}
where the conditional probability on the right-hand side can be represented as an integral:
\begin{align*}
  \P(\tau - \sigma(h_u(x,\r)) > x \mid  \sigma(h_u(x,\r)) < \tau) \hspace{-146pt} & \\
    &= \int_{h_u(x,\r)}^\infty \P(\tau - \sigma(h_u(x,\r)) > x \mid  \sigma(h_u(x,\r)) < \tau ; X(\sigma(h_u(x,\r))) = y ) \\
      &\qquad \times \dd \P(X(\sigma(h_u(x,\r))) \leq y \mid \sigma(h_u(x,\r)) < \tau) \\
    &\geq \int_{h_u(x,\r)}^\infty \P(X(t)>0; 0\leq t \leq x \mid X(0) = y) \dd \P(X(\sigma(h_u(x,\r))) \leq y \mid \sigma(h_u(x,\r)) < \tau).
\end{align*}
As the integrand is increasing in $y$, we obtain 
\begin{equation}
  \P(\tau>x; M_\tau>(1-\r)x) \geq \P(X(t)>0; 0\leq t \leq x \mid X(0) = h_u(x,\r)) \P(M_\tau> h_u(x,\r)).
  \label{eq:PisC:C:lowerbound}
\end{equation}

Rewriting the first probability on the right-hand side of \eqref{eq:PisC:C:lowerbound} yields
\begin{align*}
  \P(X(t)>0; 0\leq t \leq x \mid X(0) = h_u(x,\r)) \hspace{-162pt} & \hspace{162pt}
    = \P\left( \inf_{t\in[0,x]} \left\{X(t)-X(0)\right\} > -h_u(x,\r)\right) \\
    &\geq \P\left( \inf_{t\in[0,x]} \left\{-\r t+\sum_{i=1}^{N(t)} B_i \right\} > -g(x,\r)\right)
    = 1- \P\left( \sup_{t\in[0,x]} \left\{\r t-\sum_{i=1}^{N(t)} B_i \right\} \geq g(x,\r)\right).
\end{align*}
From Etemadi's inequality \citep[Theorem 22.5]{billingsley1995probability}, it then follows that 
\begin{align*}
  \P(X(t)>0; 0\leq t \leq x \mid X(0) = h_u(x,\r))
    &\geq 1- 3 \sup_{t\in[0,x]} \P\left(\r t-\sum_{i=1}^{N(t)} B_i \geq g(x,\r)/3\right) \\
    &\geq 1- 3 \sup_{t\in[0,x]} \P\left(\left|\r t-\sum_{i=1}^{N(t)} B_i\right| \geq g(x,\r)/3\right).
\end{align*}
The variance of $\sum_{i=1}^{N(t)} B_i$ equals $\l \E[B^2] t$ and is dominated by $2\E[B^*] t$ for all $\rho\in[0,1]$. Therefore, noting that $3\cdot 3^2\cdot 2=54$, Chebyshev's inequality implies
\begin{align*}
  \P(X(t)>0; 0\leq t \leq x \mid X(0) = h_u(x,\r))
    &\geq 1- \sup_{t\in[0,x]} \frac{54 \E[B^*] t}{g(x,\r)^2}
    = 1- \frac{54 \E[B^*]}{((1-\r)^2x)^{2p-1}} \rightarrow 1 \numberthis \label{eq:PisC:C:lowerbound2}
\end{align*}
for $x\geq x_\r^*$ and $\rhotoone$. Let $\zeta(\r):= 1- 54 \E[B^*]\left(\log\frac{1}{1-\r}\right)^{-(2p-1)k^*}$. By \eqref{eq:PisC:C:lowerbound}, \eqref{eq:PisC:C:lowerbound2} and Theorem~\ref{thm:CisNB}, one readily finds
\begin{align*}
  \sup_{x\geq x_\r^*} \P(\tau>x\mid M_\tau>(1-\r)x) 
    &\geq \sup_{x\geq x_\r^*} \P(X(t)>0; 0\leq t \leq x \mid X(0) = h_u(x,\r)) \; \frac{\P(M_\tau> h_u(x,\r))}{\P(M_\tau>(1-\r)x)} \\
    &\gtrsim \zeta(\r) \sup_{x\geq x_\r^*} \frac{\P(B> h_u(x,\r))}{\P(B>(1-\r)x)}
    \sim \zeta(\r) \left(1 + ((1-\r)^2x_\r^*)^{p-1}\right)^{-\a}, 
\end{align*}
which tends to one as $\rhotoone$. This validates expression~\eqref{eq:PandCisC} in Theorem~\ref{thm:PisC}.

\subsection{Asymptotics of unconditional first hitting time}
\label{subsec:PisCupperbound}
This section proves expression \eqref{eq:PandC}. From \eqref{eq:PandCisC} it follows that $\lim_\rhotoone \inf_{x\geq x_\r^*} \P(\tau>x \mid M_\tau>(1-\r)x) = 1$, and thus
\begin{equation}
  \lim_\rhotoone \inf_{x\geq x_\r^*} \frac{\P(\tau>x)}{\P(M_\tau>(1-\r)x)} = 1.
  \label{eq:PandCisC:lowerbound2}
\end{equation}
Proving $\lim_\rhotoone \sup_{x\geq x_\r^*} \frac{\P(\tau>x)}{\P(M_\tau>(1-\r)x)} = 1,$ however, requires far more work.

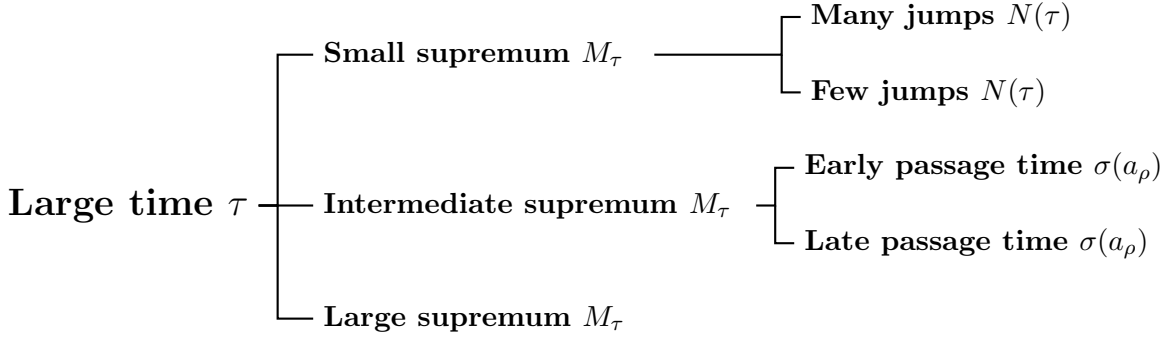
\begin{figure}[t]
\centering
 \begin{tikzpicture}[x=2.5cm, y=2cm, every node/.style={font=\bfseries}]
  \node (1) at (-.1,.25) {\Large{Large time $\tau$}};
  \node (2) at (2,1.25) {\parbox{\widthof{Intermediate supremum $M_\tau$}}{Small supremum $M_\tau$}};
  \node (3) at (2,.25) {\parbox{\widthof{Intermediate supremum $M_\tau$}}{Intermediate supremum $M_\tau$}};
  \node (4) at (2,-.5) {\parbox{\widthof{Intermediate supremum $M_\tau$}}{Large supremum $M_\tau$}};
  \node (2a) at (4.4,1.5) {\parbox{\widthof{Early passage time $\s(a)$}}{Many jumps $N(\tau)$}};
  \node (2b) at (4.4,1) {\parbox{\widthof{Early passage time $\s(a)$}}{Few jumps $N(\tau)$}};
  \node (3a) at (4.4,.5) {\parbox{\widthof{Early passage time $\s(a_\r)$}}{Early passage time $\s(a_\r)$}};
  \node (3b) at (4.4,0) {\parbox{\widthof{Early passage time $\s(a_\r)$}}{Late passage time $\s(a_\r)$}};
  \draw[thick] (1.east) --++ (.1,0) |- (2.west);
  \draw[thick] (1.east) --++ (.1,0) |- (3.west);
  \draw[thick] (1.east) --++ (.1,0) |- (4.west);
  \draw[thick] (2a.west) --++ (-.1,0) --++ (0,-.25) --++ (-.67,0);
  \draw[thick] (2b.west) --++ (-.1,0) --++ (0,.25);
  \draw[thick] (3a.west) --++ (-.1,0) --++ (0,-.25) --++ (-.1,0);
  \draw[thick] (3b.west) --++ (-.1,0) --++ (0,.25);
 \end{tikzpicture}
 \caption{Visualization of proof structure. The event of a large time $\tau$ is analysed under three scenarios, depending on the size of the supremum $M_\tau$. Two of these scenarios are again considered in more detail, where a distinction is based on the number of jumps before $\tau$ and the passage time of a high level $a_\r$.}
 \label{fig:structurePproof}
\end{figure}

As noted at the beginning of this section, the event $\{\tau>x\}$ is analysed by discriminating various scenarios. First, we specify scenarios $\{\tau>x,M_\tau \in \cdot\}$, where the supremum $M_\tau$ can be in three regions: small, intermediate and large. Then, the small and intermediate regions are shown to be negligible in Sections~\ref{subsubsec:smallC} through~\ref{subsubsec:intermediateC}.\rev{REP2COM30} Finally, the large $M_\tau$ region is shown to be asymptotically equivalent to $\P(\tau>x)$ in Section~\ref{subsubsec:largeC}. The structure of the proof is visualized in Figure~\ref{fig:structurePproof}.

We now formalize the various scenarios. Fix constants $\e_\g>0$ and $\e_\d\in\left(0,(p-\frac{1}{2})k^*-1\right)$, and define the functions $\g_\r:=\left(\log\frac{1}{1-\r}\right)^{-\e_\g}$ and $\d_\r:=\left(\log\frac{1}{1-\r}\right)^{-(1+\e_\d)}$. Similar as before, the function
$
 h_l(x,\r):= (1-\g_\r)(1-\r)x - g(x,\r),
$
where
$
 g(x,\r) = (1-\r)^{2p-1}x^p,
$
represents a lower bound for, yet is asymptotically equivalent to, $(1-\r)x$. The three regions for $M_\tau$ are now given by
\begin{align*}
  \P(\tau>x) 
    &\leq \P(\tau>x; M_\tau\leq \d_\r(1-\r)x) + \P(\tau>x; M_\tau\in (\d_\r(1-\r)x,h_l(x,\r)]) + \P(M_\tau>h_l(x,\r)) \\
    &= \mathrm{I} + \mathrm{II} + \mathrm{III}.
\end{align*}
The next paragraphs show that terms $\mathrm{I}$ and $\mathrm{II}$ both vanish faster than $\P(M_\tau>(1-\r)x)$ for $x\geq x_\r^*$ as $\rhotoone$. Contrastingly, the final paragraph shows that term $\mathrm{III}$ asymptotically behaves as $\P(M_\tau>(1-\r)x)$ in the same limiting regime.

\subsubsection[Small supremum Mtau: term I]{Small supremum $M_\tau$: term $\mathrm{I}$}
\label{subsubsec:smallC}
Term $\mathrm{I}$ is the probability of a large first hitting time $\tau$ for which the corresponding process supremum $M_\tau$ is relatively small. First, we show that the number of jumps before $\tau$ is not much higher than the expected number of jumps. Then, we show that it is highly unlikely  for a probable amount of small jumps to incur a large $\tau$.

Recall that $N(t)$ denotes the number of jumps during an interval of length $t$. In particular, $N(t)$ is Poisson distributed with mean $\l t$. Let $N_0(t)$ be the number of jumps of size at most $\d_\r(1-\r)x$ and $N_1(t)$ be the number of jumps of at least that size. Then $N_0(t)$ is Poisson distributed with mean $\l t \P(B<\d_\r(1-\r)x)$, $N_1(t)$ is Poisson distributed with mean $\l t \P(B\geq \d_\r(1-\r)x)$ and $N(t) = N_0(t) + N_1(t)$ for all $t\geq 0$. Let $B_{0,i}$ be i.i.d. random variables with c.d.f.\ $\P(B_{0,i} \leq y) = \P(B \leq y \mid B\leq \d_\r(1-\r)x)$.\rev{REP2COM32}

If $\tau>x$, all jumps before time $x$ had a cumulative size of at least $x$; that is, $\sum_{i=0}^{N(x)} B_i>x$ is a necessary condition for $\tau>x$. Furthermore, the inequality $M_\tau \geq B_\tau$ is trivial. Let $\l^*=(1+\eta_\r)\l$, where $\eta_\r=(1-\r)/2$. Note that this implies $\l^*\E[B] \in [0,1)$ whenever $\rho\in [0,1)$. Now
\begin{align*}
  \P(\tau>x, M_\tau \leq \d_\r(1-\r)x) \hspace{-27pt} & \hspace{27pt}
      = \P\left(\tau>x, \sum_{i=0}^{N(x)} B_i>x, M_\tau \leq \d_\r(1-\r)x \right) \\
    &\leq \P\left(\sum_{i=0}^{N(x)} B_i>x, \bigvee_{i=0}^{N(x)} B_i \leq \d_\r(1-\r)x\right) 
      = \P\left(\sum_{i=0}^{N_0(x)} B_{0,i}>x, N_1(x) = 0 \right) \\
    &\leq \P\left(\sum_{i=0}^{N_0(x)} B_{0,i}>x \right)
      \leq \P(N_0(x)\geq \l^* x) + \P\left(\sum_{i=0}^{\l^* x} B_{0,i}>x\right) \\
    &= \P(N(x)\geq \l^* x) + \P\left(\sum_{i=0}^{\l^* x} B_i>x \bigg| \bigvee_{i=0}^{\l^* x} B_i \leq \d_\r(1-\r)x\right) \\
    &= \mathrm{Ia} + \mathrm{Ib}.
\end{align*}
Term $\mathrm{Ia}$ corresponds to a system where the number of jumps greatly exceeds its expectation. Term $\mathrm{Ib}$ indicates a likely number of jumps, none of which has a size exceeding $\d_\r(1-\r)x$. 

\subsubsection[Many jumps: term Ia]{Many jumps: term $\mathrm{Ia}$}
From Markov's inequality, one can see that for all $s\geq 0$ we have
\begin{equation*}
  \mathrm{Ia} = \P(e^{sN(x)}\geq e^{s\l^* x}) \leq e^{-s\l^* x} \E\left[e^{sN(x)}\right] = \exp[-\l x((1+\eta_\r)s - e^{s}+1)].
\end{equation*}
Taking the infimum over all $s\geq 0$ gives
\begin{equation*}
  \mathrm{Ia} \leq \exp[-\l x\sup_{s\geq 0}((1+\eta_\r) s - e^{s}+1)] = \exp\left[-\l x\left((1+\eta_\r)\log(1+\eta_\r) - \eta_\r\right)\right].
\end{equation*}
The bound $\log(1+\eta_\r)\geq \frac{2\eta_\r}{2+\eta_\r}$ for $\eta_\r>0$ then yields
$\mathrm{Ia} \leq \exp\left[-\frac{\eta_\r^2}{2+\eta_\r}\l x\right].$
Dividing by $\P(M_\tau>(1-\r)x)$, taking the supremum, applying Theorem~\ref{thm:CisNB} and using Potter's Theorem with $\nu>0$ gives
\begin{align*}
  \sup_{x\geq x_\r^*} \frac{\mathrm{Ia}}{\P(M_\tau>(1-\r)x)} 
    &\lesssim \sup_{x\geq x_\r^*} C \frac{1-\r}{\r\P(B>(1-\r)x)}\exp\left[-\frac{\eta_\r^2}{2+\eta_\r}\l x\right] \\
    &\leq \sup_{x\geq x_\r^*} C \frac{1-\r}{\r}\exp\left[(\a+\nu)\log((1-\r)x)-\frac{\eta_\r^2}{2+\eta_\r}\l x\right] \\
    &\leq C \frac{1-\r}{\r}\exp\left[(\a+\nu)\log\frac{1}{1-\r}-\frac{\l}{10}\log^{k^*}\frac{1}{1-\r} + o\left(\log\frac{1}{1-\r}\right)\right] \\
    &\rightarrow 0, \numberthis \label{eq:PisC:Ia}
\end{align*}
as $\rhotoone$.

\subsubsection[Few jumps: term Ib]{Few jumps: term $\mathrm{Ib}$}
Now consider term $\mathrm{Ib}$. The corresponding event is a large $\tau$, caused by a probable amount of small jumps. The following theorem by \citet{prokhorov1959extremal} is used to show that this scenario is unlikely as $\rho$ tends to 1.
\begin{theorem}[\citealp{prokhorov1959extremal}, Theorem 1] \label{thm:prokhorov}
 Suppose that $\xi_i,i=1,\ldots,n$ are independent, zero-mean random variables such that there exists a constant $c$ for which $|\xi_i|\leq c$ for $i=1,\ldots,n$, and $\sum_{i=1}^n \Var\{\xi_i\} < \infty$. Then
 \begin{equation}
  \P\left(\sum_{i=1}^n \xi_i > y\right) \leq \exp\left[-\frac{y}{2c}\arcsinh\frac{yc}{2\sum_{i=1}^n \Var\{\xi_i\}} \right].
  \label{eq:prokhorovoriginal}
 \end{equation}
\end{theorem}
Using the bound $\arcsinh(z) = \log (z+\sqrt{1+z^2}) \geq \log(2z)$, Prokhorov's inequality implies
\begin{equation*}
 \P\left(\sum_{i=1}^n \xi_i > y\right) \leq \left(\frac{cy}{\sum_{i=1}^n \Var\{\xi_i\}}\right)^{-\frac{y}{2c}}.
 \tag{\ref{eq:prokhorovoriginal}*}
 \label{eq:prokhorov}
\end{equation*}
Define $Y_i := B_i-\E[B]$. Then
\begin{align*}
 \mathrm{Ib} &= \P\left(\sum_{i=0}^{\l^* x} Y_i> (1-\l^*\E[B])x - \E[B] \bigg| \bigvee_{i=0}^{\l^* x} Y_i \leq \d_\r(1-\r)x - \E[B] \right).
\end{align*}\rev{REP1COM25}\rev{REP2COM33}
Let $\sigma^2_{B_0}$ be the variance of $B$ provided $B<\d_\r(1-\r)x$. Then $\sigma^2_{B_0} \leq \sigma^2_B$ and hence, using~\eqref{eq:prokhorov},
\begin{align*}
  \mathrm{Ib} &\leq \left( \frac{\d_\r(1-\r)x - \E[B]}{\l^* x+1}\frac{(1-\l^*\E[B])x - \E[B]}{\sigma^2_B}\right)^{-\frac{(1-\l^*\E[B])x-\E[B]}{2\d_\r(1-\r)x-2\E[B]}} \\
    &= \exp\left[-(1+\phi_\r^{(1)}(x))\frac{1-\l^*\E[B]}{2\d_\r(1-\r)}\log\left(\frac{1-\phi_\r^{(2)}(x)}{\l^* \sigma^2_B}(1-\l^*\E[B])(1-\r)\d_\r x\right)\right],
\end{align*}
where the real-valued functions $\phi^{(i)}_\r(x)$ are defined as
\begin{equation*}
 \phi_\r^{(1)}(x):=\frac{1-\frac{\E[B]}{(1-\l^*\E[B])x}}{1-\frac{\E[B]}{\d_\r(1-\r)x}}-1,
 \qquad
 \phi_\r^{(2)}(x):= 1-\frac{\left(1 - \frac{\E[B]}{\d_\r(1-\r)x}\right)\left(1-\frac{\E[B]}{(1-\l^*\E[B])x}\right)}{1+\frac{1}{\l^* x}},
\end{equation*}
and satisfy $\phi_\r^{(i)}(x)\rightarrow 0$ as $\rhotoone$ for $x\geq x_\r^*$. Additionally, the functions $\phi_\r^{(i)}$ are non-negative and non-increasing for $\rho$ sufficiently close to one. These properties imply that the inequality
\begin{align*}
  \mathrm{Ib} &\leq \exp\left[-\frac{1-\l^*\E[B]}{2\d_\r(1-\r)}\log\left(\frac{1-\phi_\r^{(2)}(x)}{\l^* \sigma^2_B}(1-\l^*\E[B])(1-\r)\d_\r x\right)\right],
\end{align*}
holds for $\rho$ sufficiently close to one and $x\geq x_\r^*$. Substitution of $\lambda^*=(1+\eta_\r)\l = \frac{3-\r}{2}\l$ subsequently gives
\begin{align*}
  \mathrm{Ib} 
   &\leq \exp\left[-\frac{1}{4\d_\r}\log\left(\frac{1-\phi_\r^{(2)}(x)}{3\l \sigma^2_B}(1-\r)^2\d_\r x\right)\right].
\end{align*}
Dividing the upper bound above by $\P(M_\tau>(1-\r)x) \sim \frac{\r}{1-\r}\P(B>(1-\r)x)$ and applying Potter's Theorem with $\nu>0$ yields
\begin{align*}
  \frac{\mathrm{Ib}}{\P(M_\tau>(1-\r)x)} 
    &\lesssim C \frac{1-\r}{\r} \exp\bigg[(\a+\nu)\log((1-\r)x) - \frac{1}{4\d_\r}\log\left(\frac{1-\phi_\r^{(2)}(x)}{3 \l \sigma^2_B}(1-\r)^2 \d_\r x\right) \bigg] \\
    &= C \frac{1-\r}{\r} \exp\bigg[\left(\a+\nu-\frac{1}{4\d_\r}\right)\log((1-\r)x) - \frac{1}{4\d_\r}\log\left(\frac{1-\phi_\r^{(2)}(x)}{3 \l \sigma^2_B}(1-\r)\d_\r\right) \bigg].
\end{align*}
The supremum over $x\geq x_\r^*$ is attained in $x=x_\r^*$ for $\r$ sufficiently close to one. That is,
\begin{align*}
  \sup_{x\geq x_\r^*} \frac{\mathrm{Ib}}{\P(M_\tau>(1-\r)x)} &\lesssim C \frac{1-\r}{\r} \exp\bigg[\left(\a+\nu-\frac{1}{4\d_\r}\right)\log\left(\frac{1}{1-\r}\log^{k^*}\frac{1}{1-\r}\right) \\
      &\qquad - \frac{1}{4\d_\r}\log\left(\frac{1-\phi_\r^{(2)}(x_\r^*)}{3 \l \sigma^2_B}(1-\r)\d_\r\right) \bigg] \\
    &= C \frac{1-\r}{\r} \exp\bigg[(\a+\nu)\log\left(\frac{1}{1-\r}\log^{k^*}\frac{1}{1-\r}\right) \\
      &\qquad - \frac{1}{4}\log^{1+\e_\d}\left(\frac{1}{1-\r}\right)\log\left(\frac{1-\phi_\r^{(2)}(x_\r^*)}{3 \l \sigma^2_B} \log^{k^*-1-\e_\d}\frac{1}{1-\r}\right) \bigg] \\
    &\rightarrow 0 \numberthis \label{eq:PisC:Ib}
\end{align*}
as $\rhotoone$. Together, \eqref{eq:PisC:Ia} and \eqref{eq:PisC:Ib} assure that term $\mathrm{I}$ is dominated by $\P(M_\tau>(1-\r)x)$.

\subsubsection[Intermediate supremum Mtau: term II]{Intermediate supremum $M_\tau$: term $\mathrm{II}$}
\label{subsubsec:intermediateC}
Term $\mathrm{II}$ corresponds to the event of a large $\tau$ that experiences an intermediate supremum $M_\tau$\rev{REP2COM34}. Write
\begin{align*}
  \sup_{x\geq x_\r^*} \frac{\mathrm{II}}{\P(M_\tau>(1-\r)x)} 
    &\leq \sup_{x\geq x_\r^*} \frac{\P(M_\tau > \d_\r(1-\r)x)}{\P(M_\tau > (1-\r)x)} \\
      &\qquad \times \sup_{x\geq x_\r^*} \P(\tau>x; M_\tau \leq h_l(x,\r) \mid M_\tau > \d_\r(1-\r)x) \\
    &\sim \d_\r^{-\a} \sup_{x\geq x_\r^*} \P(\tau>x; M_\tau \leq h_l(x,\r) \mid M_\tau > \d_\r(1-\r)x). \numberthis \label{eq:PisC:II}
\end{align*}
Set $\kappa_\r := \left(\log\frac{1}{1-\r}\right)^{-\e_\kappa}$ for some $\e_\kappa\geq \e_\g$, implying $\g_\r-\kappa_\r>0$. By considering the time $\s(a_\r)$ when the process $X(t)$ first exceeds level $a_\r:=\d_\r(1-\r)x$, we can partition~\eqref{eq:PisC:II} into two events:
\begin{align*}
  \P(\tau>x; M_\tau \leq h_l(x,\r) \mid M_\tau > \d_\r(1-\r)x) \hspace{-148pt} & \hspace{148pt} 
    = \P(\tau>x; \s(a_\r)\leq \kappa_\r x; M_\tau \leq h_l(x,\r) \mid \s(a_\r) < \tau) \\
    &\qquad + \P(\tau>x; \s(a_\r)>\kappa_\r x; M_\tau \leq h_l(x,\r) \mid \s(a_\r) < \tau) \\
    &\leq \P(\tau>(1-\kappa_\r)x; M_\tau \leq h_l(x,\r) \mid \s(a_\r) < \tau) + \P(\s(a_\r)>\kappa_\r x \mid \s(a_\r) < \tau) \\
    &= \mathrm{IIa} + \mathrm{IIb}.
\end{align*}
Term $\mathrm{IIa}$ is associated with sample paths that experiences an intermediate supremum and that may already hit zero after time $(1-\kappa_\r)x$. Term $\mathrm{IIb}$ corresponds to a sample path where the process does not exceed level $a_\r$ before time $\kappa_\r x$, provided that it will hit level $a_\r$ before it hits zero. 

\subsubsection[Early passage time: term IIa]{Early passage time: term $\mathrm{IIa}$}
\label{subsec:earlypassage}
Term $\mathrm{IIa}$ is analysed along the lines of Section~\ref{subsec:PisClowerbound}:
\begin{align*}
  \mathrm{IIa} 
    &= \int_0^{h_l(x,\r)} \P(X(t)>0; 0\leq t \leq (1-\kappa_\r) x; M_\tau \leq h_l(x,\r) \mid \s(a_\r) < \tau; X(0) = y) \\
      &\qquad \times \dd \P(X(0) \leq y \mid \s(a_\r) < \tau) \\
    &\leq \int_0^{h_l(x,\r)} \P(X(t)>0; 0\leq t \leq (1-\kappa_\r) x \mid \s(a_\r) < \tau; X(0) = y)  \dd \P(X(0) \leq y \mid \s(a_\r) < \tau) \\
    &\leq \P(X(t)>0; 0\leq t \leq (1-\kappa_\r) x \mid X(0) = h_l(x,\r)) \\
    &\leq \P(X((1-\kappa_\r) x)-X(0)>-h_l(x,\r)). 
\end{align*}
Here, the second inequality holds as the integrand is increasing in $y$, and $a_\r\leq h_l(x,\r)$ for $x\geq x_\r^*$ and $\r$ sufficiently close to one.

Define $A^\r_0:=0$ and $A^\r_i:=\inf\{t\geq 0: N(\sum_{j=0}^{i-1}A^\r_j+t)\geq i\}$ for all $i\geq 1$. Then the $A^\r_i$ are i.i.d.\ exponentially distributed random variables with mean $1/\l$ and $\sum_{i=1}^{N(t)} A^\r_i \leq t$ for all $t\geq 0$. We drop the superscript $\r$ for notational convenience. Now,
\begin{align*}
  \mathrm{IIa} &\leq \P\left(-(1-\kappa_\r)x + \sum_{i=1}^{N((1-\kappa_\r) x)} B_i >-h_l(x,\r)\right) \\
    &= \P\left(-\r(1-\kappa_\r)x + \sum_{i=1}^{N((1-\kappa_\r) x)} B_i > (\g_\r-\kappa_\r)(1-\r)x + g(x,\r)\right) \\
    &\leq \P\left(\sum_{i=1}^{N((1-\kappa_\r) x)} [B_i-\r A_i] > (\g_\r-\kappa_\r)(1-\r)x + g(x,\r)\right).
\end{align*}
Fix $q\in\left(\max\left\{2,\frac{(1+\e_\d)\a}{(p-\frac{1}{2})k^*}\right\}, \a\right)$. By Chebyshev's inequality for general moments and Theorem 5.1 in Chapter 1 of \citet{gut1988stopped}, there exists some constant $C_q$ such that
\begin{align*}
  \mathrm{IIa} &\leq \frac{\E\left[\left(\sum_{i=1}^{N((1-\kappa_\r) x)} [B_i-\r A_i]\right)^q\right]}{\left((\g_\r-\kappa_\r)(1-\r)x + g(x,\r)\right)^q}
    \leq \frac{C_q \E[|B_1-\r A_1|^q] \E[N((1-\kappa_\r) x)^{q/2}]}{g(x,\r)^q} \\
    &\leq \frac{C_q \E[|B_1-\r A_1|^q] \E[N((1-\kappa_\r) x)]^{q/2}}{g(x,\r)^q},
\end{align*}
where the last derivation is justified by H{\"o}lder's inequality. Subsequently, one may show from Jensen's inequality that $\E[|B-A|^q]\leq 2^{q-1}(\E[|A|^q]+\E[|B|^q])$ and therefore
\begin{align*}
  \mathrm{IIa} 
    &\leq \frac{C_q 2^{q-1} (\E[B_1^q]+\E[A_1^q]) (\l (1-\kappa_\r) x)^{q/2}}{(1-\r)^{(2p-1)q}x^{pq}}
    \leq C \left((1-\r)^2x\right)^{-(p-\frac{1}{2})q}.
\end{align*}
for some constant $C$. By choice of $p,q$ and $\e_\d$, we conclude that
\begin{equation}
  \d_\r^{-\a} \sup_{x\geq x_\r^*} \mathrm{IIa}
    \leq C \left(\log\frac{1}{1-\r}\right)^{(1+\e_\d)\a-(p-\frac{1}{2})k^*q}
    \rightarrow 0
  \label{eq:IIa}
\end{equation}
as $\rhotoone$.

\subsubsection[Late passage time: term IIb]{Late passage time: term $\mathrm{IIb}$}
Term $\mathrm{IIb}$ is analysed with the following crucial lemma, and is proven in Section~\ref{sec:firstpassagetime}:
\begin{lemma}\label{lem:LateJump}
 Suppose $\P(B>x) = L(x)x^{-\a}$ for some $\a>2,\a\neq 3$ and $L(x)$ slowly varying. Define $a_\r^*:= k^* \mu (\a-1)\frac{1}{1-\r}\log\frac{1}{1-\r}$ for some $k^*>2$. Then for any fixed $y>0$,
\begin{equation}
 \sup_{a\geq a_\r^*} \E[\s(a)\mid \s(a)<\tau; X_0=y] = \O\left(\frac{1}{1-\r}\right)
 \label{eq:FixedJumpTime}
\end{equation}
as $\rhotoone$. Similarly, without conditioning on the value of $X_0$,
\begin{equation}
 \sup_{a\geq a_\r^*} \E[\s(a)\mid \s(a)<\tau] = \O\left(\frac{1}{1-\r}\right)
 \label{eq:RandomJumpTime}
\end{equation}
as $\rhotoone$.\rev{REP2COM2}
\end{lemma}
Applying Markov's inequality and sequentially Lemma~\ref{lem:LateJump} to term $\mathrm{IIb}$ yields, as $\rhotoone$,
\begin{align*}
  \d_\r^{-\a} \sup_{x\geq x_\r^*} \mathrm{IIb} 
    &\leq \sup_{x\geq x_\r^*} \frac{\E[\s(\d_\r(1-\r)x)\mid \s(\d_\r(1-\r)x)<\tau]}{\d_\r^\a \kappa_\r x} \\
    &= \O\left(\frac{1}{1-\r}\right) \frac{(1-\r)^2}{\left(\log\frac{1}{1-\r}\right)^{k^*-(1+\e_\d)\a-\e_\kappa}}
    \rightarrow 0. \numberthis\label{eq:IIb}
\end{align*}

\subsubsection[Large supremum Mtau: term III]{Large supremum $M_{\tau}$: term $\mathrm{III}$}
\label{subsubsec:largeC}
Finally, we show that the probability of a large time $\tau$ is asymptotically equivalent to term $\mathrm{III}$. Using Theorem~\ref{thm:CisNB}, it directly follows that
\begin{align*}
  \sup_{x\geq x_\r^*} \frac{\P(M_\tau>h_l(x,\r))}{\P(M_\tau>(1-\r)x)} 
    &\lesssim \sup_{x\geq x_\r^*} \frac{\P(B>h_l(x,\r))}{\P(B>(1-\r)x)} \sim \sup_{x\geq x_\r^*} \left(\frac{(1-\g_\r)(1-\r)x-(1-\r)^{2p-1}x^p}{(1-\r)x}\right)^{-\a} \\
    &= \left(1-\left(\log\frac{1}{1-\r}\right)^{-\e_\g}-\left(\log\frac{1}{1-\r}\right)^{-(1-p)k^*}\right)^{-\a}
    \rightarrow 1
\end{align*}\rev{REP2COM31}
as $\rhotoone$. This completes the proof of~\eqref{eq:PandC}.

\section[Asymptotics of the conditional expectation of the passage time of level a]{Asymptotics of the conditional expectation of the passage time of level $a$}
\label{sec:firstpassagetime}
This section is dedicated to the proof of Lemma~\ref{lem:LateJump}, which regards the expected first passage time of level $a$, $\s(a)$, provided that level $a$ is reached before level $0$: $\s(a)<\tau$. In particular, we consider high levels $a\geq a_\r^*:=k^*\mu(\a-1)\frac{1}{1-\r}\log\frac{1}{1-\r}$ for any $k^*>2$. The lemma considers two different scenarios. In the first scenario, we condition on the initial value $X(0)=y$. In the second scenario, the initial value $X(0)$ is a random variable with the same distribution as a general jump size $B$. The analysis for this latter scenario is based on the following decomposition:
\begin{equation}
 \E[\s(a)\mid\s(a)<\tau] = \int_0^a \E[\s(a)\mid\s(a)<\tau; X(0)=y] \dd \P(B\leq y).\rev{REP2COM19}
 \label{eq:firstpassagetime}
\end{equation}
That is, we condition the former expectation to the initial value of the process and integrate over all possible initial values. A distinction is made between a ``small'' and a ``large'' random initial value; a precise definition of which is given at the end of these introductory paragraphs. The first scenario, where the initial value is fixed, is implicit in the analysis of a small random initial value, and its proof is concluded at the end of Section~\ref{subsec:smallinitialvalue}.

The derivation of results in this section relies heavily on the theory of spectrally one-sided L{\'e}vy processes and $q$-scale functions, e.g.\ as documented by \citet{kyprianou2014introductory}. Our interest in $q$-scale functions $W_\r^{(q)}$ originates from the close connection between the all-time supremum $M_\infty$ and the $0$-scale function $W_\r(x):=W_\r^{(0)}(x)$. Of particular importance is the relation\rev{REP2COM13}
\begin{equation}
 \P(M_\infty < x) = (1-\r)W_\r(x),
 \label{eq:Wtoprob}
\end{equation}
which can be derived from Corollary IX.3.4 in \citet{asmussen2003applied} \citep[e.g.\ as shown in][]{bekker2009advances}. Prior to stating the definition of $W_\r^{(q)}$\rev{REP1COM17}, we define the Laplace exponent $\psi(\l):=\frac{1}{t}\log \E(e^{-\l X(t)})$ 
and the right inverse $\varphi(q):=\sup\{\l\geq 0:\psi(\l)=q\}$. Now, for every $q\geq 0$ the $q$-scale function $W_\r^{(q)}(x): \R\rightarrow [0,\infty)$ corresponding to the spectrally positive L{\'e}vy process $X(t)$ is defined on $x<0$ as $W_\r^{(q)}(x)=0$, and on $x\geq 0$ by its Laplace transform:
\begin{equation}
 \int_0^\infty e^{-\beta x}W_\r^{(q)}(x) \dd x = \frac{1}{\psi(\beta)-q} \text{ for } \beta >\varphi(q).
\end{equation}
Additionally, \citeauthor{kyprianou2014introductory} gives a representation of $W_\r^{(q)}(x)$ in terms of $W_\r(x)$ in his equation (8.29):
\begin{equation}
 W_\r^{(q)}(x) = \sum_{k\geq 0} q^k W_\r^{(k+1)\circledast}(x),
 \label{eq:convolutionWq}
\end{equation}
where the function $f^{1\circledast}(x)$ is identical to $f(x)$ and $f^{k\circledast}(x):=\int_0^x f^{(k-1)\circledast}(x-y)f(y) \dd y$ denotes the $k$-fold convolution of $f$ with itself.\rev{REP1COM18}

An alternative representation of $W_\r(x)$ is provided by expression (8.22) in \citet{kyprianou2014introductory}\rev{REP1COM19}\rev{REP2COM20}, stating that there are a measure $n_\r(\cdot)$ on the space of excursions of $X(t)$ from its previous minimum $\min\{X(s):0\leq s\leq t\}$ and a random variable $\bar{\xi}_\r$ associated with the height of an excursion, such that for all $b>x\geq 0$ we have\rev{REP2COM21}
\begin{equation}
 W_\r(x) = W_\r(b)\exp\left(-\int_x^b n_\r(\bar{\xi}_\r>t) \dd t\right).
 \label{eq:Wexponentform}
\end{equation}
This representation will provide a useful property for the all-time supremum p.d.f.\ $f_{M_\infty}(x) := \frac{\dd}{\dd y}\P(M_\infty<y)\bigg|_{y=x}$. Using the Pollaczek-Khintchine formula (cf. equation \ref{eq:WisS}), we write
\begin{align*}
  f_{M_\infty}(x) &= \sum_{n=1}^\infty (1-\r)\r^n \frac{\dd}{\dd y} \P(B_1^* + \ldots + B_n^* < y)\bigg|_{y=x}
\end{align*}
for $x>0$. One may show by induction that $\frac{\dd}{\dd y}\P(B_1^* + \ldots + B_n^* < y)$ is defined everywhere and is bounded by $1/\E[B]$ for all $n\geq 1$. As such, $f_{M_\infty}(x)$ is properly defined and bounded for all $x>0$.\rev{REP1COM20} Additionally, \eqref{eq:Wexponentform} implies that
\begin{equation}
 \frac{f_{M_\infty}(x)}{\P(M_\infty<x)} = \frac{\dd}{\dd y}\log W_\r(y)\bigg|_{y=x} = n_\r(\bar{\xi}_\r>x)
 \label{eq:logW}
\end{equation}
is non-increasing in $x$.

For the remainder of this section, the subscripts $\r$ for $W_\r(x)$ and $W_\r^{(q)}(x)$ are discarded. We also introduce the short-hand notations $\E_y[\cdot]$ and $\P_y(\cdot)$ for the conditional expectation $\E[\cdot\mid X(0)=y]$ and conditional probability $\P(\cdot\mid X(0)=y)$, respectively.

Let $Z^{(q)}(x) := 1+q\int_0^x W^{(q)}(y) \dd y$. From \eqref{eq:convolutionWq} and the spectrally \textit{positive} L{\'e}vy process interpretation of Theorem 8.1 in \citet{kyprianou2014introductory}, it follows that
\begin{align*}
 \E_y[\s(a) \,\ind\{\s(a)<\tau\} ] \hspace{-46pt}&\hspace{46pt}= -\frac{d}{dq} \E_y[e^{-q\s(a)} \ind\{\s(a)<\tau\} ]\bigg|_{q=0} \\
  &= -\frac{d}{dq} Z^{(q)}(a-y) + \frac{W^{(q)}(a-y)}{W^{(q)}(a)} \frac{d}{dq} Z^{(q)}(a) \\
   &\qquad + Z^{(q)}(a) \frac{W^{(q)}(a)\frac{d}{dq} W^{(q)}(a-y) - W^{(q)}(a-y) \frac{d}{dq} W^{(q)}(a)}{(W^{(q)}(a))^2} \bigg|_{q=0} \\
  &= -\int_0^{a-y} W(t) \dd t + \frac{W(a-y)}{W(a)} \int_0^a W(t) \dd t + \frac{W(a)W^{2\circledast}(a-y) - W(a-y)W^{2\circledast}(a)}{(W(a))^2} \\
  &= \frac{W(a-y)}{W(a)} \frac{\int_0^a (W(a)-W(a-t))W(t) \dd t}{W(a)} - \frac{\int_0^{a-y}(W(a)-W(a-y-t))W(t)\dd t}{W(a)}.
\end{align*}
Now, from (8.12) in \citet{kyprianou2014introductory}\rev{REP1COM21} one may deduce $\P_y(\s(a)<\tau) = \frac{W(a)-W(a-y)}{W(a)}$, which gives a representation of the conditional expectation in terms of scale functions:
\begin{multline*}
 \E_y[\s(a) \mid \s(a) < \tau] = \frac{W(a-y)}{W(a)} \frac{\int_0^a (W(a)-W(a-t))W(t) \dd t}{W(a)-W(a-y)} \\
  - \frac{\int_0^{a-y}(W(a)-W(a-y-t))W(t)\dd t}{W(a)-W(a-y)}.
\end{multline*}
Substitute \eqref{eq:Wtoprob} into the above expression to obtain
\begin{align*}
  \E_y[\s(a) \mid \s(a) < \tau] \hspace{-36pt}&\hspace{36pt}= \frac{\P(M_\infty < a-y)}{\P(M_\infty < a)}\frac{\int_0^a \P(M_\infty \in [a-t,a))\P(M_\infty < t) \dd t}{(1-\r)\P(M_\infty \in [a-y,a))} \\
    &\hspace{36pt}\qquad - \frac{\int_0^{a-y}\P(M_\infty \in [a-y-t,a))\P(M_\infty < t) \dd t}{(1-\r)\P(M_\infty \in [a-y,a))} \\
  &\leq \frac{\int_0^a \P(M_\infty \in [a-t,a))\P(M_\infty < t) \dd t - \int_0^{a-y}\P(M_\infty \in [a-y-t,a))\P(M_\infty < t) \dd t}{(1-\r)\P(M_\infty \in [a-y,a))} \\
  &=: \frac{K_{num}(y,a)}{K_{denom}(y,a)} \numberthis \label{eq:expectationtau}
\end{align*}
The analysis of this expression depends on the initial value $y$. We distinguish two categories of initial values: small and large values. Fix $d$ such that $0<d<1-\frac{2}{k^*}<1$. Small values are of size at most $d\cdot a$, all other values are large values.

\subsection{Small random initial value or fixed initial value}
\label{subsec:smallinitialvalue}
This section considers the process from a small initial value $y$, i.e.\ $y\leq da$. For any $y$-differentiable function $G(y,a)$, it is known that $G(y,a) = G(0,a) + \int_0^y \frac{d}{ds} G(s,a)\big|_{s=z} \dd z$. This is now used to obtain an alternative representation of $K_{num}(y,a)$.

Let $M_\infty^{(i)}, i=1,2$ be independent copies of $M_\infty$. Taking the derivative of $K_{num}(s,a)$ with respect to $s$ yields
\begin{align*}
 \frac{d}{ds} K_{num}(s,a) &= \P(M_\infty^{(2)} < a)\P(M_\infty^{(1)} < a-s) - \int_0^{a-s} \P(M_\infty^{(1)} < t) \dd \P(M_\infty^{(2)} < a-s-t) \\
  &= \P(M_\infty^{(2)} < a)\P(M_\infty^{(1)} < a-s) - \P(M_\infty^{(1)} + M_\infty^{(2)} < a-s) \\
  &= \P(M_\infty^{(1)} < a-s) - \P(M_\infty^{(1)} + M_\infty^{(2)} < a-s) - \P(M_\infty^{(2)} \geq a)\P(M_\infty^{(1)} < a-s) \\
  &= \P(M_\infty^{(1)} + M_\infty^{(2)} \geq a-s; M_\infty^{(1)} < a-s) \\
   &\qquad - \P(M_\infty^{(2)} \geq a-s)\P(M_\infty^{(1)} < a-s) + \P(M_\infty^{(2)} \in [a-s,a))\P(M_\infty^{(1)} < a-s) \\
  &= \P(M_\infty^{(1)} + M_\infty^{(2)} \geq a-s; M_\infty^{(1)} < a-s; M_\infty^{(2)} < a-s) \\
   &\qquad + \P(M_\infty^{(2)} \in [a-s,a))\P(M_\infty^{(1)} < a-s),
\end{align*}
so that $K_{num}(0,a)=0$ implies\rev{REP1COM22}\rev{REP2COM22}
\begin{align*}
 \E_y[\s(a) \mid \s(a) < \tau] &\leq \frac{K_{num}(y,a)}{K_{denom}(y,a)} \\
  &= \frac{\int_0^y \P(M_\infty^{(1)} + M_\infty^{(2)} \geq a-z; M_\infty^{(1)} < a-z; M_\infty^{(2)} < a-z) \dd z}{(1-\r)\P(M_\infty \in [a-y,a))} \\
   &\qquad + \frac{\int_0^y \P(M_\infty^{(2)} \in [a-z,a))\P(M_\infty^{(1)} < a-z) \dd z}{(1-\r)\P(M_\infty \in [a-y,a))} \\
  &\leq \frac{\int_0^y \P(M_\infty^{(1)} + M_\infty^{(2)} \geq a-z; M_\infty^{(1)} < a-z; M_\infty^{(2)} < a-z) \dd z}{(1-\r)\P(M_\infty \in [a-y,a))} + \frac{y}{1-\r}
\end{align*}
By symmetry, we have
\begin{align*}
  \P(M_\infty^{(1)} + M_\infty^{(2)} \geq u; M_\infty^{(1)} < u; M_\infty^{(2)} < u)
    &\leq 2\P(M_\infty^{(1)} + M_\infty^{(2)} \geq u; u/2 \leq M_\infty^{(1)} < u) \\
    &\leq 2\P(M_\infty \in [u/2,u))
\end{align*}
and hence
\begin{align*}
  \E_y[\s(a) \mid \s(a) < \tau] &\leq \frac{2 \int_0^y \P\left(M_\infty \in \left[\frac{a-z}{2},a-z\right)\right) \dd z}{(1-\r)\P(M_\infty \in [a-y,a))} + \frac{y}{1-\r} \\
    &\leq \frac{2y}{1-\r}\left(1+\frac{\P\left(M_\infty \in \left[\frac{a-y}{2},a \right)\right)}{\P(M_\infty \in [a-y,a))}\right) \numberthis \label{eq:Mratio}
\end{align*}

Both local probabilities can be represented as a sum of local probabilities over an interval with fixed length. Subsequently, Theorem~\ref{thm:WisNB*} is applied to bound the above ratio. Fix $y_{min}>0$ and first consider \eqref{eq:Mratio} for $y_{min}\leq y\leq da$. For $S:=y_{min}/2$, we have
\begin{align*}
 \frac{\P\left(M_\infty \in \left[\frac{a-y}{2},\frac{a}{2}\right)\right)}{\P\left(M_\infty \in \left[a-y,a\right)\right)}
  &\leq \frac{\sum_{i=0}^{\lceil \frac{y}{2S}-1\rceil} \P\left(M_\infty \in \left[\frac{a-y}{2}+iS,\frac{a-y}{2}+(i+1)S\right)\right)}{\sum_{i=0}^{\lfloor \frac{y}{S}-1\rfloor} \P\left(M_\infty \in \left[a-y+iS,a-y+(i+1)S\right)\right)}.
\end{align*}
We would now like to utilise Theorem~\ref{thm:WisNB*}. To this end, consider $x_\r$ as defined by Theorem~\ref{thm:WisNB*} with parameter $(1-d)k^*/2>1$. Then for all $y\leq da$, we have $\frac{a-y}{2}\geq \frac{1-d}{2} a_\r^*=x_\r$. Hence, we observe that there exists a non-increasing function $\phi_\r(\cdot)\downarrow 0$ for which the inequalities\rev{REP2COM23}
\begin{equation}
 1-\phi_\r\left(\frac{a-y}{2}\right) \leq \frac{\P\left(M_\infty \in \left[\frac{a-y}{2}+iS,\frac{a-y}{2}+(i+1)S\right)\right)}{\frac{\r}{1-\r}\P\left(B^* \in \left[\frac{a-y}{2}+iS,\frac{a-y}{2}+(i+1)S\right)\right)} \leq 1+\phi_\r\left(\frac{a-y}{2}\right)
\end{equation}
both hold for all $y\leq da$ and $i\geq 0$. From $a-y \geq a-da \geq \frac{a}{k^*}$ one may subsequently conclude that the ratio of interest is bounded:
\begin{align*}
  \frac{\P\left(M_\infty \in \left[\frac{a-y}{2},\frac{a}{2}\right)\right)}{\P\left(M_\infty \in \left[a-y,a\right)\right)}
  &\leq \frac{1+\phi_\r\left(\frac{a}{2k^*}\right)}{1-\phi_\r\left(\frac{a}{k^*}\right)} \frac{\sum_{i=0}^{\lceil \frac{y}{2S}-1\rceil} \P\left(B^* \in \left[\frac{a-y}{2}+iS,\frac{a-y}{2}+(i+1)S\right)\right)}{\sum_{i=0}^{\lfloor \frac{y}{S}-1\rfloor} \P\left(B^* \in \left[a-y+iS,a-y+(i+1)S\right)\right)} \\
  &\leq \frac{1+\phi_\r\left(\frac{a}{2k^*}\right)}{1-\phi_\r\left(\frac{a}{k^*}\right)} \frac{\frac{y}{2S}+1}{\frac{y}{S}-1} \frac{\P\left(B > \frac{a-y}{2}\right)}{\P\left(B>a \right)}
  \sim \frac{1+\frac{2S}{y}}{2-\frac{2S}{y}} \left(\frac{a-y}{2a}\right)^{-\a}
  \leq 2 (2k^*)^{\a}.
\end{align*}

Second, consider \eqref{eq:Mratio} for $0<y<y_{min}$. Relation~\eqref{eq:logW} implies
\begin{align*}
 \frac{\P\left(M_\infty \in \left[\frac{a-y}{2},\frac{a}{2}\right)\right)}{\P\left(M_\infty \in \left[a-y,a\right)\right)}
  &\leq 
  \frac{y\sup_{z\in (0,y)} \frac{f_{M_\infty}\left(\frac{a-z}{2}\right)}{\P\left(M_\infty < \frac{a-z}{2}\right)}\P\left(M_\infty < \frac{a-z}{2}\right)}{2y \inf_{z\in (0,y)} \frac{f_{M_\infty}\left(a-z\right)}{\P(M_\infty < a-z)}\P(M_\infty < a-z)} \leq \frac{\frac{f_{M_\infty}\left(\frac{a-y}{2}\right)}{\P\left(M_\infty < \frac{a-y}{2}\right)}\P\left(M_\infty < \frac{a}{2}\right)}{2 \frac{f_{M_\infty} (a)}{\P(M_\infty < a)}\P(M_\infty < a-y)} \\
  &= \frac{f_{M_\infty}\left(\frac{a-y}{2}\right)}{2 f_{M_\infty} (a)} \frac{\P\left(M_\infty < \frac{a}{2}\right)\P(M_\infty < a)}{\P\left(M_\infty < \frac{a-y}{2}\right)\P(M_\infty < a-y)}\sim \frac{f_{M_\infty}\left(\frac{a-y}{2}\right)}{2 f_{M_\infty} (a)}
\end{align*}
as $a\rightarrow\infty$. We conclude that
\begin{equation}
 \E_y[\s(a) \mid \s(a) < \tau] \lesssim C\frac{y}{1-\r}\left(1 + \frac{f_{M_\infty}\left(\frac{a-y}{2}\right)}{f_{M_\infty} (a)} \ind\{y\leq y_{min}\} \right). \label{eq:firstpassagetimesmallB}
\end{equation}
The above relation explicitly shows the dependence of the asymptotic upper bound on $y$. This dependence is crucial in the analysis of the second part of the lemma, where we will integrate the upper bound over $\P(B<y)$. However, before addressing large initial values it should be noted that \eqref{eq:firstpassagetimesmallB} also proves the first part of the lemma. There, $y$ is fixed and the lemma follows directly after choosing $0<y_{min}<y$.\rev{REP2COM24} 

\subsection{Large random initial value}
Complementary to the previous section, we now consider \eqref{eq:expectationtau} for large initial values, i.e.\ $da \leq y < a$. Let $M_\infty^*$ be a random variable with the excess distribution of $M_\infty$ as its c.d.f.\ , that is, $\frac{d}{dx} \P(M_\infty^* < x) = \P(M_\infty \geq t)/\E[M_\infty]$. Using $\P(M_\infty < t)=1-\P(M_\infty \geq t)$ and $\int_0^a \P(M_\infty \in [a-t,a)) \dd t = \E[M_\infty \ind\{M_\infty < a\}]$, we find
\begin{align*}
  K_{num}(y,a)
    &= \int_0^a \P(M_\infty \in [a-t,a))\P(M_\infty < t) \dd t - \int_0^{a-y}\P(M_\infty \in [a-y-t,a))\P(M_\infty < t) \dd t \\
    &= \E[M_\infty\ind\{M_\infty < a\}] - \E[M_\infty] \int_0^a \P(M_\infty \in [a-t,a)) \dd \P(M_\infty^* < t) \\
      &\qquad - \E[M_\infty\ind\{M_\infty < a-y\}] + \E[M_\infty] \int_0^{a-y} \P(M_\infty \in [a-y-t,a)) \dd \P(M_\infty^* < t) \\
    &\leq \E[M_\infty\ind\{M_\infty \in [a-y,a)\}] + \E[M_\infty] \P(M_\infty \in [a-y-M_\infty^*,a); M_\infty^* < a-y) \\
    &\leq a\P(M_\infty \in [a-y,a)) + \E[M_\infty].
\end{align*}
It therefore follows that
\[
  \E_y[\s(a) \mid \s(a) < \tau] - \frac{a}{1-\r} \leq \frac{\E[M_\infty]}{(1-\r)\P(M_\infty \in [a-y,a))} \leq  \frac{\E[M_\infty]}{(1-\r)\P(M_\infty \in [(1-d)a,a))},
\]
where $\E[M_\infty] = \frac{\r}{1-\r}\frac{\E[B^2]}{2\E[B]}$. Similar to the analysis of small initial values, Theorem~\ref{thm:WisNB*} invokes
\begin{align*}
  \E_y[\s(a) \mid \s(a) < \tau]
    &\lesssim \frac{a}{1-\r} + \frac{C}{(1-\r)da\P\left(B > a \right)}. \numberthis \label{eq:firstpassagetimelargeB}
\end{align*}
This completes the analysis of the conditional expectation for large initial values.

\subsection{Synthesis of small and large random initial value}
From equation \eqref{eq:firstpassagetime}, 
\eqref{eq:firstpassagetimesmallB} 
and \eqref{eq:firstpassagetimelargeB} one can deduce that
\begin{align*}
 \sup_{a\geq a_\r^*} \E[\s(a)\mid \s(a)<\tau]\hspace{-48pt}& \\
  &\leq \sup_{a\geq a_\r^*} \int_0^{da} \E_y[\s(a)\mid\s(a)<\tau] \dd \P(B < y)  + \sup_{a\geq a_\r^*} \int_{da}^a \E_y[\s(a)\mid\s(a)<\tau] \dd \P(B < y) \\
  &\lesssim \frac{C}{1-\r} \sup_{a\geq a_\r^*} \int_0^{da} y \dd \P(B < y) + \frac{C}{1-\r} \sup_{a\geq a_\r^*} \int_0^{y_{min}} y\cdot \frac{f_{M_\infty}\left(\frac{a-y}{2}\right)}{f_{M_\infty} (a)} \dd\P(B<y) \\
   &\qquad + \sup_{a\geq a_\r^*} \P(B \geq da) \sup_{y\in [da,a)} \E_y[\s(a)\mid\s(a)<\tau] \\
  &\lesssim \frac{C\E[B]}{1-\r} + \frac{C y_{min}}{1-\r} \sup_{a\geq a_\r^*} \int_0^{y_{min}} \frac{f_{M_\infty}\left(\frac{a-y}{2}\right)}{f_{M_\infty} (a)} \dd \P(B < y) \\\rev{REP2COM25}
   &\qquad + \sup_{a\geq a_\r^*} \frac{a}{1-\r}\P(B \geq da) + \sup_{a\geq a_\r^*} \frac{C}{(1-\r) a} \frac{\P(B \geq da)}{\P\left(B \geq a\right)}. \numberthis \label{eq:anyJumpI}
\end{align*}
The third term is dominated by its Markov's bound $\frac{\E[B]}{(1-\r)d}$. Also, the integral in the second term is ultimately bounded by a constant. This follows from the fact that $\frac{f_{M_\infty}(x)}{\P(M_\infty\leq x)}$ is non-increasing and application of Theorem~\ref{thm:WisNB*} as before:
\begin{align*}
  \int_0^{y_{min}} \frac{f_{M_\infty}\left(\frac{a-y}{2}\right)}{f_{M_\infty} (a)} \dd \P(B<y) \hspace{-26pt} & \hspace{26pt}
    \leq \frac{\P(M_\infty\leq \frac{a}{2})}{f_{M_\infty} (a)} \int_0^{y_{min}} \frac{f_{M_\infty}\left(\frac{a-y}{2}\right)}{\P(M_\infty\leq \frac{a-y}{2})} \dd \P(B<y) \\
    &\leq \P(B<y_{min}) \frac{\P(M_\infty\leq \frac{a}{2})}{\P(M_\infty\leq a)} \frac{\P(M_\infty\leq a)}{f_{M_\infty} (a)} \frac{f_{M_\infty}\left(\frac{a-y_{min}}{2}\right)}{\P(M_\infty\leq \frac{a-y_{min}}{2})} \\
    &= C \frac{\P(M_\infty\leq \frac{a}{2})}{\P(M_\infty\leq a)} \inf_{y\in(0,y_{min})} \frac{\P(M_\infty\leq a+y)}{f_{M_\infty} (a+y)} \inf_{y\in(0,y_{min})} \frac{f_{M_\infty}\left(\frac{a+y-2y_{min}}{2}\right)}{\P(M_\infty\leq \frac{a+y-2y_{min}}{2})} \\
    &\leq C \frac{\P(M_\infty\leq \frac{a}{2})}{\P(M_\infty\leq a)} \frac{\P(M_\infty\leq a+y_{min})}{\P(M_\infty\leq \frac{a-y_{min}}{2})} \frac{\inf_{y\in(0,y_{min})} f_{M_\infty}\left(\frac{a+y-2y_{min}}{2}\right)}{\sup_{y\in(0,y_{min})} f_{M_\infty} (a+y)} \\
    &\lesssim C \frac{\int_0^{y_{min}} f_{M_\infty}\left(\frac{a+y-2y_{min}}{2}\right) \dd y}{\int_0^{y_{min}} f_{M_\infty} (a+y) \dd y}  = C \frac{\P\left(M_\infty \in \left(\frac{a-2y_{min}}{2},\frac{a-y_{min}}{2}\right)\right)}{\P(M_\infty\in(a,a+y_{min}))} \\
    &\lesssim C \frac{\P\left(B > \frac{a-2y_{min}}{2}\right)}{\P(B > a+y_{min})} \sim C \left(1-\frac{3y_{min}}{a+y_{min}}\right)^{-\a}
\end{align*}
as $a\rightarrow \infty$. Substituting this into~\eqref{eq:anyJumpI} gives\rev{REP2COM22}
\begin{equation}
  \sup_{a\geq a_\r^*} \E[\s(a)\mid \s(a)<\tau] \lesssim \frac{C}{1-\r} + \frac{C y_{min}}{1-\r} \sup_{a\geq a_\r^*}  \left(1-\frac{3y_{min}}{a+y_{min}}\right)^{-\a} + \sup_{a\geq a_\r^*} \frac{C}{(1-\r)a}d^{-\a}.
\end{equation}
Since all suprema are obtained in $a=a_\r^*$ as $\rhotoone$
, the above expressions can be written in terms of $1/(1-\r)$:
\begin{align*}
  \sup_{a\geq a_\r^*} \E[\s(a)\mid \s(a)<\tau]
    &\lesssim \frac{C}{1-\r} + \frac{C}{\log\frac{1}{1-\r}} 
    = \O\left(\frac{1}{1-\r}\right), \hspace{1cm} \rhotoone.
\end{align*}

\section{Tightness of bounds -- proofs}
\label{sec:tightness}
This section presents the proof of Lemma~\ref{lem:Wdensity}, Corollary~\ref{cor:lowerboundxrho} and Lemma~\ref{lem:lowerboundxrhostar}, respectively.

\subsection{Local Kingman heavy traffic approximation}
Complete monotonicity of $f_{M_\infty}(\cdot)$ follows from Corollary 3.2 in \citet{keilson1978exponential}. As $f_{M_\infty}(\cdot)$ is non-increasing, it follows that the random variable $V$ with c.d.f.\ $F_V(0):=1, F_V(x):= 1-\frac{f_{M_\infty}(x)}{f_{M_\infty}(0+)}, x>0,$ is well-defined. Relation \eqref{eq:Wdensity} is now derived by analysing the Laplace-Stieltjes transform of $V$.
Let $\tilde{M}_\infty(\cdot)$ and $\tilde{B^*}(\cdot)$ denote the Laplace-Stieltjes transforms of $M_\infty$ and $B^*$, respectively. One the one hand, we have \citep[equation (7.9)]{adan2002queueing}
\begin{equation}
  \int_{0+}^\infty e^{-st} f_{M_\infty}(t) \dd t = \tilde{M}_\infty(s) - \P(M_\infty = 0)= \frac{1-\r}{1-\r \tilde{B^*}(s)} - (1-\r) = \frac{\r(1-\r)\tilde{B^*}(s)}{1-\r\tilde{B^*}(s)}.  \label{eq:Wdensity1}
\end{equation}
On the other hand, integration by parts yields
\begin{align*}
  \int_{0+}^\infty e^{-st} f_{M_\infty}(t) \dd t
     &= \frac{1}{s} f_{M_\infty}(0+) + \frac{1}{s} \int_{0+}^\infty e^{-st} \dd f_{M_\infty}(t). \numberthis \label{eq:Wdensity2}
\end{align*}
Combining \eqref{eq:Wdensity1} and \eqref{eq:Wdensity2} gives
\[
  \E[e^{sV}] = -\int_0^\infty e^{-st} \dd \frac{f_{M_\infty}(t)}{f_{M_\infty}(0+)} = 1-\frac{\r(1-\r)s\tilde{B^*}(s)}{f_{M_\infty}(0+)(1-\r\tilde{B^*}(s))}= 1-\frac{\E[B]s\tilde{B^*}(s)}{1-\r\tilde{B^*}(s)},
\]
since $f_{M_\infty}(0+) = (1-\r)\lambda$ (cf.\ relation \ref{eq:WisS}). Now
\begin{align*}
  \E[e^{(1-\r)sV}] 
    &= 1-\frac{\E[B](1-\r)s\tilde{B^*}((1-\r)s)}{1-\r\tilde{B^*}((1-\r)s)} 
     = 1-\frac{\E[B](1-\r)s\tilde{B^*}((1-\r)s)}{1-\r\left(1-\E[B^*](1-\r)s + o(1-\r)\right)} \\
    &\rightarrow 1-\frac{\E[B]s}{1+\E[B^*]s}
\end{align*}
as $\rhotoone$. 
Inverting this expression and applying the Continuity Theorem gives
\begin{equation}
  \P((1-\r)V \leq x) \rightarrow 1-\frac{\E[B]}{\E[B^*]} e^{-\frac{x}{\E[B^*]}},
\end{equation}
provided $\E[B^*]\geq \E[B]$. Under this assumption, the lemma statement follows from the definition of $F_V(x)$. The proof is therefore concluded once we verify that all completely monotone densities $f_B(\cdot)$ satisfy $\E[B^*]\geq \E[B]$.

Bernstein's theorem \citep{bernstein1929fonctions} 
states that any completely monotone function can be represented as mixture of exponential functions. In particular, there exists a non-decreasing function $\mu(\cdot)$ such that
\begin{equation}
  f_B(x) = \int_0^\infty e^{-tx} \dd \mu(t).
\end{equation}
From this representation, one may derive $1 = \int_0^\infty \frac{1}{t} \dd \mu(t)$, $\E[B] = \int_0^\infty \frac{1}{t^2} \dd \mu(t)$ and $\E[B^2] = \int_0^\infty \frac{2}{t^3} \dd \mu(t)$.
A straightforward computation yields 
\[
 \E[B^2] - 2\E[B]^2
 = \int_0^\infty \int_0^\infty \frac{1}{st}\left(\frac{1}{s} - \frac{1}{t}\right)^2 \dd \mu(s) \dd \mu(t) \geq 0.
\]
The claimed property follows from $\E[B^*]-\E[B] = (\E[B^2]-2\E[B]^2)/(2\E[B])\geq 0$.

\subsection[Lower bound of the function xrho]{Lower bound of the function $x_\r$}
Since the p.d.f.\ of both $M_\infty$ and $B^*$ are well-defined and non-increasing, one can see that
\begin{align*}
  \frac{\P(M_\infty \in \frac{y}{1-\r} + \D)}{\frac{\r}{1-\r}\P(B^* \in \frac{y}{1-\r} + \D)}
    &= \frac{\int_{y/(1-\rho)}^{y/(1-\rho)+T} f_{M_\infty}(t) \dd t}{\frac{\r}{1-\r} \int_{y/(1-\rho)}^{y/(1-\rho)+T} \P(B>t)/\E[B] \dd t}
    \geq \frac{f_{M_\infty}\left(\frac{y}{1-\rho}+T\right)}{\frac{\l}{1-\r} \P\left(B>\frac{y}{1-\rho}\right)} \\
    &\geq \frac{\frac{1}{1-\r} f_{M_\infty}\left(\frac{y+T}{1-\rho}\right)}{\frac{\l}{(1-\r)^2} \P\left(B>\frac{y}{1-\rho}\right)}.
\end{align*}
Fix $0<\nu<\a-2$. According to Potter's Theorem there exists a constant $C>0$ such that $\P(B>x)\leq C x^{-\a+\nu}$ for $x$ sufficiently large. Hence, by Lemma~\ref{lem:Wdensity},
\begin{align*}
  \lim_\rhotoone \frac{\P\left(M_\infty \in \frac{y}{1-\r} + \D\right)}{\frac{\r}{1-\r}\P\left(B^* \in \frac{y}{1-\r} + \D\right)}
    &\geq \lim_\rhotoone \frac{\frac{1}{1-\r} f_{M_\infty}\left(\frac{y+T}{1-\rho}\right)}{\l C (1-\r)^{\a-2-\nu} y^{-\a+\nu}}=\infty.
\end{align*}

\subsection[Lower bound of the function xrhostar]{Lower bound of the function $x_\r^*$}
Theorem~1 in \citet{abate1995limits} yields
\begin{equation}
  \lim_\rhotoone \frac{1}{1-\r}\frac{\E[B^2]}{\E[B]^2}\P\left(\tau>\frac{\E[B^2]t}{\E[B](1-\r)^2}\right) =t^{-1/2} \sqrt{\frac{2}{\pi}} e^{-t/2} - 2 \bar{\Phi}(t^{1/2})=: f_{\mathrm{R}}(t).
\end{equation}
Thus, as $F_B$ is regularly varying with index $-\a<-2$,
\begin{align*}
  \lim_{\rhotoone} \frac{\P(\tau > y_\r^*)}{\frac{\r}{1-\r} \P(B > (1-\r)y_\r^*)}
    &= \frac{\E[B]^2}{\E[B^2]} \lim_{\rhotoone} \frac{\frac{1}{1-\r}\frac{\E[B^2]}{\E[B]^2} \P\left(\tau > \frac{y}{(1-\r)^2}\right)}{ \frac{\r}{(1-\r)^2} \P\left(B > \frac{y}{1-\r}\right)} \\
    &= \frac{\E[B]^2}{\E[B^2]} f_R\left(\frac{\E[B]^2}{\E[B^2]} y\right) \lim_{\rhotoone} \frac{1}{ \frac{\r}{(1-\r)^2} \P\left(B > \frac{y}{1-\r}\right)}=\infty.
\end{align*}

\bibliographystyle{apalike}
\DeclareRobustCommand{\NLprefix}[3]{#3}

\end{document}